\documentclass[11pt]{amsart}
\usepackage{amsmath,amssymb,verbatim}

\newtheorem{thm}{Theorem}[section]
\newtheorem{lemma}[thm]{Lemma}
\newtheorem{prop}[thm]{Proposition}
\newtheorem{cor}[thm]{Corollary} 

\DeclareMathOperator{\charp}{char}
\DeclareMathOperator\Aut{Aut}
\DeclareMathOperator\Gal{Gal}
\DeclareMathOperator\PSL{PSL_2}
\DeclareMathOperator\PGL{PGL_2}
\DeclareMathOperator\SL{SL_2}
\DeclareMathOperator\GL{GL_2}
\DeclareMathOperator\GLn{GL}
\DeclareMathOperator\ind{ind}
\DeclareMathOperator\Hom{Hom}
\DeclareMathOperator\PGU{PGU_3}
\DeclareMathOperator\PGammaL{P{\Gamma}L_2}
\newcommand{\Cone}{{\mathcal C}}
\newcommand{\Ctwo}{{\mathcal D}}
\newcommand{\Cthree}{{\mathcal B}}
\DeclareMathOperator\orb{orb}

\newcommand\textmatrix[1]{\bigl(\begin{smallmatrix}#1\end{smallmatrix}\bigr)}

\renewcommand{\theenumi}{{\roman{enumi}}}

\newcommand\Line{{\mathbb P}^1}
\newcommand\F{{\mathbb{F}}}
\newcommand\bZ{{\mathbb{Z}}}
\newcommand\C{{\mathbb{C}}}

\renewcommand{\bar}[1]{#1\llap{$\overline{\phantom{\rm#1}}$}}
\newcommand{\col}{\,{:}\,}

\title{Polynomials with $\operatorname{PSL}(2)$ Monodromy}

\author{Robert M. Guralnick}
\address{Department of Mathematics, University of Southern California,
Los Angeles, CA 90089-2532, USA}
\email{guralnic@usc.edu}

\author{Michael E. Zieve}
\address{
Center for Communications Research, 805 Bunn Drive,
 Princeton, NJ 08540-1966, USA}
\email{zieve@math.rutgers.edu}
\urladdr{http://www.math.rutgers.edu/$\sim$zieve/}

\date{\today}

\thanks{We thank Joel Rosenberg and the referee for catching some
typographical errors.
The first author was partially supported by NSF grant DMS 0653873}

\begin{document}

\begin{abstract}
Let $k$ be a field of characteristic $p>0$, let
$q$ be a power of $p$, and let $u$ be transcendental over $k$.
We determine all polynomials
$f\in k[X]\setminus k[X^p]$ of degree $q(q-1)/2$ for which the
Galois group of $f(X)-u$ over $k(u)$ has a transitive normal
subgroup isomorphic to $\PSL(q)$, subject to a certain
ramification hypothesis.  As a consequence, we describe all
polynomials $f\in k[X]$ such that $\deg(f)$ is not a power of $p$
and $f$ is functionally indecomposable over $k$ but $f$ decomposes over
an extension of $k$.  Moreover, except for one ramification configuration
(which is handled in a companion paper with Rosenberg), we describe all
indecomposable
polynomials $f\in k[X]$ such that $\deg(f)$ is not a power of $p$
and $f$ is exceptional, in the sense that $X-Y$ is the only absolutely
irreducible factor of $f(X)-f(Y)$ which lies in $k[X,Y]$.
It is known that, when $k$ is finite, a polynomial $f$ is exceptional
if and only if it induces a bijection on infinitely many finite extensions
of $k$.
\end{abstract}

\maketitle


\section{Introduction}

Let $\Cone$ and $\Ctwo$ be smooth, projective, geometrically irreducible curves
over a field $k$ of characteristic $p\ge 0$, and let $f\colon \Cone\to \Ctwo$ be a
separable morphism over $k$ of degree $d>1$.  Much information about the
map $f$ is encoded in its \emph{monodromy groups},
which are defined as follows.  Let $k(\Cone)/k(\Ctwo)$ be the separable field
extension corresponding to $f$, and let $E$ denote its Galois closure.
The \emph{arithmetic monodromy group} of $f$ is the group
$A:=\Gal(E/k(\Ctwo))$.
Letting $\ell$ denote the algebraic closure of $k$ in $E$, the
\emph{geometric monodromy group} of $f$ is $G:=\Gal(E/\ell(\Ctwo))$.

A fundamental problem is to determine the possibilities for the
monodromy groups and the ramification of such maps $f$, where $\Ctwo$ is
fixed and $\Cone$ (and $f$) varies.  Riemann solved this problem in case
$k=\C$, and moreover he determined how many such maps $f$ have a
specified monodromy group and specified branch points and inertia groups
(although it remains unknown how to write down equations for these $f$).
Riemann's result can be generalized to the case that $k$ is any
algebraically closed field of characteristic zero; but the problem becomes
much more difficult
for other fields $k$.  In case $k$ is an algebraically closed field of
characteristic $p>0$, the best result to date was proved by
Raynaud \cite{ray} and Harbater {\cite{har}, and describes the
geometric monodromy groups of maps $f$ whose branch points are contained
in a specified set.  However, the problem of determining the
possible inertia groups (not to mention the higher ramification groups)
is wide open.

In this paper we explicitly determine all maps $f$ (in positive
characteristic) having certain monodromy groups, subject to a constraint
on the ramification.
The specific situation we consider arises in questions about two
special types of maps $f$.  In many applications, one is interested
in maps $f$ satisfying some additional hypotheses; in practice, one
often finds that natural hypotheses on $f$ imply severe restrictions
on the monodromy groups.  Classifying the groups satisfying such
conditions often requires deep group theoretic results (by contrast,
no difficult group theory was involved in the proofs of Raynaud and
Harbater).  We need some notation to describe our situation.
The map $f$ induces a transitive permutation representation of $A$
on a set $\Delta$ of size $d$.  Let $A_1$ be a point stabilizer in this
representation, and note that $G$ is also transitive on $\Delta$.

In this paper we study maps $f$ having some of the following conditions;
some of these conditions are geometric and some arithmetic, and all of them
have been studied for over a century.  In our list we include translations
of the conditions into properties of the monodromy groups.

\begin{enumerate}
\item $f$ is (arithmetically) indecomposable (i.e., $f$ is not a nontrivial
composition of maps defined over $k$).  This is equivalent to $A_1$ being a
maximal subgroup of $A$.
\item $f$ is geometrically indecomposable.  This is equivalent
to $A_1 \cap G$ being maximal in $G$.
\item  $f$ is totally ramified at some point. This says an inertia group
in the Galois closure of $f$ acts transitively on $\Delta$.
\item $\Cone$ has genus $g$. This can be translated to a property of $G$ and
the higher ramification groups, via the Riemann-Hurwitz formula and
Hilbert's different formula.
\item The fiber product $\{(c,d) \in \Cone \times \Ctwo : f(c)=f(d)\}$ has no
geometric components defined over $k$ except the diagonal.
This says that $A$ and $G$ have no
common orbits on $\Delta \times \Delta$ besides the diagonal.
\end{enumerate}

In case (v) we say $f$ is an \emph{exceptional map}.  These have been
studied extensively (starting with Dickson's 1896 thesis \cite{Dith}, and
subsequently by Schur, Carlitz, Davenport, Lewis, Bombieri, Fried,
Cohen, and others).  They are particularly
interesting when $k$ is finite.   In this case,  $f$ is exceptional
if and only if $f$ is bijective on $k'$-rational points for all extensions
$k'/k$ of degree relatively prime to some positive integer $m$
(which can always be taken to be $|A\col G|$).  Indeed, as long as 
the cardinality of $k$ is sufficiently large compared to the degree of $f$
and the genus of $\Cone$,
exceptionality is equivalent to either injectivity or surjectivity of
the map induced by $f$ on $k$-rational points.  Moreover, for finite $k$,
the composition of two maps $\Cone\to \Cthree$ and $\Cthree\to \Ctwo$ is exceptional if and
only if both maps are exceptional; thus, it suffices to classify the
exceptional maps that are arithmetically indecomposable.
For proofs of these results, see \cite[p.~2]{GTZ} and \cite[p.~185]{FGS}.
For partial analogues over infinite fields, see \cite{GS2,LMZ,Mac}.

We will be especially interested in the case that $f$ is a polynomial;
here $\Cone$ has genus $g=0$, and $f$ is totally ramified over
a rational point of $\Ctwo$ (namely, the point at infinity).

In this paper (except for a few cases handled in 
\cite{GRZ,GRZ2}), we classify the polynomials $f(X)\in k[X]$ such that
$\deg(f)$ is not a power of $\charp(k)$ and at least one of the
following holds:

\begin{enumerate}
\item[(1)] $f$ is arithmetically but not geometrically indecomposable;
\item[(2)] $f$ is exceptional and arithmetically indecomposable; or
\item[(3)] $A$ has a transitive normal subgroup isomorphic to $\PSL(q)$.
\end{enumerate}

It was shown in \cite{GS} that there are significant restrictions on
the monodromy groups of an arithmetically indecomposable polynomial.
In this paper we will study the groups that can occur for polynomials $f$
which satisfy either (1) or (2).  We will see that these usually
give rise to condition (3).  It would be of great value to have
a classification of all indecomposable $f$ for which $G$ is neither
alternating nor symmetric; from \cite{GS} we know that the situation in
case (3) is one of the main sources of such polynomials $f$.
In characteristic $0$, all such $f$ are known \cite{Fe,Fe2,Mumon};
in that case there are no polynomials satisfying (1), and the
polynomials satisfying (2) have been classified \cite{Fr}.

We now consider (1) in more detail.  The problem here is to find indecomposable
polynomials over $k$ which decompose over a bigger field.
There are many examples of such polynomials in the classical
family of additive polynomials $\sum \alpha_i X^{p^i}$; further
examples occur in the related family of subadditive polynomials,
where we say $S(X)\in k[X]$ is subadditive if there is a
positive integer
$m$ and an additive polynomial $L$ such that $L(X)^m=S(X^m)$.
Up to composition with linears, these were the only examples
known before 1993.  Work of Guralnick and Saxl \cite{GS,GS2}
showed that there are severe restrictions on the degree of
any polynomial satisfying (1).  We extend and refine their result as follows,
and in particular we determine all such polynomials whose degree
is not a power of the characteristic; these include some variants
of a degree-$21$ example found by M\"uller, as well as new examples
of degree $55$.

\begin{thm}
\label{ind}
Let $k$ be a field of characteristic $p$.  If $f\in k[X]$
is indecomposable over $k$ but decomposes over some
extension of $k$, then one of the following holds:
\begin{enumerate}
\item $\deg(f)=p^e$ with $e\geq 2$;
\item $\deg(f)=21$ and $p=7$;
\item $\deg(f)=55$ and $p=11$.
\end{enumerate}
For $k$ of characteristic $p\in\{7,11\}$, there exist such $f$
of degree not a power of $p$ if and only if $k$ contains
nonsquares; moreover, all such $f$ are described in
Theorem~\ref{indthm}.
\end{thm}

We note that a new family of examples of degree $p^e$ was
found recently by Beals and Zieve \cite{BZ}, and we expect
that these examples (and the additive and subadditive examples)
will comprise all examples of degree $p^e$.

Now consider (2).  The classical examples of exceptional
polynomials are the additive
and subadditive polynomials discussed above (which are exceptional
precisely when they have no nonzero root in $k$), the multiplicative
polynomials $X^d$ (which are exceptional when $k$ contains no $d$-th
roots of unity besides $1$), and the Dickson polynomials $D_d(X,\alpha)$.
Here for $\alpha\in k$ the Dickson polynomial is defined by the equation
$D_d(Y+\alpha/Y,\alpha)=Y^d+(\alpha/Y)^d$, and its exceptionality criteria are similar
to those of $D_d(X,0)=X^d$.
All of these examples occurred in Dickson's 1896 thesis \cite{Dith}, and no
further examples were found for almost a century.  In fact, the theme
of most work in the century following Dickson's thesis was to show
that compositions of the known exceptional polynomials (including
linear polynomials) comprised all exceptional polynomials in some
situations.  Klyachko \cite{Kl} proved this for polynomials whose
degree is either equal to or coprime to $p$.
Cohen \cite{Co2p} and
Wan \cite{W2p} proved the same result for degree $2p$.
The following result of Fried, Guralnick and Saxl \cite{FGS,GS2}
provides a vast generalization of these results:

\begin{thm}
\label{fgsthm}
Let $k$ be a field
of characteristic $p$, and let $f\in k[X]$ be indecomposable and
exceptional of degree $d>1$.  Then the geometric monodromy group $G$
of $f$ satisfies one of the following conditions:
\begin{enumerate}
\item $G$ is cyclic or dihedral of odd prime degree $d\neq p$.
\item $d=p^e$ and $G=\F_p^e\rtimes G_1$, where\/ $\F_p^e$ acts on itself
by translations and $G_1\le\GLn(\F_p^e)$.
\item $p\in\{2,3\}$ and $d=p^e(p^e-1)/2$ with $e>1$ odd,
 and\/ $\PSL(p^e)$ is a transitive normal subgroup of $G$.
\end{enumerate}
\end{thm}

In particular, the degree of an indecomposable exceptional
polynomial is either (i) a prime distinct from $p$, or (ii) a power of $p$,
or (iii) $p^e(p^e-1)/2$ with $e>1$ odd and $p\in\{2,3\}$.
Any polynomial in (i) is (up to composition with linears) a Dickson
polynomial $D_d(X,\alpha)$ with $\alpha\in k$; see \cite[Appendix]{MuSchur} or
\cite{Kl}.  Case (ii) includes the additive polynomials (where $G_1=1$) and
the subadditive polynomials (where $G_1$ is cyclic).
In joint work
with M\"uller \cite{GM,GMZ}, we have found families of examples in
which $G_1$ is dihedral, and we suspect that no further examples exist in
case (ii).  This is based on the following reasoning: let $E$ denote
the Galois closure of $k(x)/k(f(x))$ (with $x$ transcendental over $k$),
and let $F$ denote the
subfield of $E$ fixed by $\F_p^e$.  We show in \cite{GMZ} that, in any
further example of (ii), the genus $g$ of $F$ would satisfy $g>1$
(whereas all known examples have $g=0$).  But then $G_1$ is a group of
automorphisms of $F$ whose order is large compared to $g$, and there
are not many possibilities for such a field $F$ \cite{GZ2}.  We hope
to complete the analysis of case (ii) in a subsequent paper.  The
present paper addresses case (iii), which does not include any classical
examples.

Case (iii) was studied intensively in the two years
following \cite{FGS}, resulting in examples with $k=\F_p$ for each
odd $e>1$ and either $p=2$ (M\"uller \cite{Mu}, Cohen--Matthews \cite{CM})
or $p=3$ (Lenstra--Zieve \cite{LZ}).
In the present paper we analyze this case in detail: we identify
all possibilities for the ramification in $k(x)/k(f(x))$,
and for all but one such possibility we determine all the corresponding
exceptional polynomials (cf.\ Theorems \ref{char3thm}
and \ref{char2thm}).  This leads to new
examples of indecomposable exceptional polynomials, which are twists
of the examples found in \cite{CM,LZ,Mu}.
In a companion paper with Rosenberg \cite{GRZ}, we
complete the analysis of case (iii) by analyzing the final
ramification possibility (which yields a new family of exceptional
polynomials).  A simplified version of Theorem~\ref{char3thm} is as
follows (and the shape of Theorem~\ref{char2thm} is similar):
\begin{thm}
\label{char3thmintro}
Let $k$ be a field of characteristic $3$, and let $q=3^e$ with $e>1$ odd.
The following are equivalent:
\begin{enumerate}
\item there exists an indecomposable exceptional polynomial $f\in k[X]$ of
degree $q(q-1)/2$;
\item $k\cap\F_q=\F_3$ and $k$ contains non-square elements.
\end{enumerate}
Moreover, these polynomials $f$ are precisely the following (up to
composition on both sides with linear polynomials in $k[X]$):
\[
X(X^{2n}-\alpha)^{(q+1)/(4n)}\left(\frac{(X^{2n}-\alpha)^{(q-1)/2}+\alpha^{(q-1)/2}}{X^{2n}}
  \right)^{(q+1)/(2n)},
\]
where $n$ divides $(q+1)/4$ and the image of $\alpha\in k^*$ in the quotient
group $k^*/(k^*)^{2n}$ has even order.
\end{thm}

For both of the above problems---finding all indecomposable 
polynomials $f(X)\in k[X]$ of degree not a power of $\charp(k)$
which either decompose over a larger field or are
exceptional---we use a similar approach.  Write $\bar{k}$ for an
algebraic closure of $k$, and let $x$ and $u$ be transcendental
over $k$.  Our general strategy
is to first translate the desired properties of the polynomial
into properties of the monodromy groups $G=\Gal(f(X)-u,\bar{k}(u))$
and $A=\Gal(f(X)-u,k(u))$, then find all group-theoretic
configurations satisfying these properties, and finally, for each
group-theoretic possibility, find all corresponding polynomials.
In our cases, a translation to group theory was done in~\cite{FGS},
and in that paper and~\cite{GS,GS2} a restricted list of plausible
pairs $(G,A)$ were given.  However, these papers did not
use the condition that $k(x)/k(f(x))$ is an extension of fields
of genus zero; via the Riemann-Hurwitz genus formula and Hilbert's different
formula, this condition leads to restrictions on the possible
ramification in the extension.  We apply this to each of the
pairs $(G,A)$ allowed for our problems by~\cite{GS2},
producing a list of all possibilities for the ramification.
The next step is to determine the possibilities for the Galois
closure $E$ of $\bar{k}(x)/\bar{k}(f(x))$; once this is
done, we compute the group $\Aut_{\overline{k}}(E)$, find all
of its subgroups which are isomorphic to $G$, and for each such
subgroup we compute the invariant subfield $E^G$ and then
compute the corresponding polynomials.  This gives all
polynomials over $\bar{k}$ having the desired
group theoretic setup {\it geometrically}; the final step is
to determine which of these polynomials are defined over $k$
and solve our original problems.

The hardest step in our work is the determination of $E$.
The data we are given for this is a group $G$ of automorphisms
of $E$, together with knowledge of the ramification in
$E/E^G$ (and the fact that $E^G$ has genus zero).
In our case, it turns out that $E$
has the shape $\bar{k}(v,w)$, where $v^{p^e}-v=w^n$
and $n$ is coprime to $p$.  So we must prove that this field
is determined by its ramification over a certain subextension;
before stating the result, we give a simple lemma describing the
ramification in the relevant subextension.

\begin{lemma}
\label{lem}
Let $k$ be a field containing\/ $\F_q$, let $n>1$ be
coprime to $q$, and let $r>0$.  For any $\gamma\in k^*$, let $v$ and $w$ be 
transcendental over $k$ such that $v^q-v=\gamma w^n$;
then the extension $k(v,w)/k(w^r)$ is Galois if and only if
$r/\!\gcd(n,r)$ divides $q-1$ and $k$ contains a primitive $r$-th root
of unity.  Moreover, for any such $r$,
the ramification is as follows (where $E=k(v,w)$ and $t=w^r$):
\begin{align}
&\text{all ramification in $E/k(t)$ occurs over two places of $k(t)$:}
  \notag\\
&\text{the finite prime 0, over which the ramification index
   is $r$;}\notag\\
&\text{and the infinite place, which is totally ramified (index
  $qr$), and} \tag{{$\dagger$}}\\
&\text{over which the sequence of ramification groups has the shape}\notag\\
&\text{ \qquad } I_0\gneqq I_1=\dots=I_n\gneqq I_{n+1}=1.\notag
\end{align}
\end{lemma}

\begin{thm}
\label{unq}
Let $k$ be a perfect field containing\/ $\F_{p^e}$,
let $t$ be transcendental over $k$, and let $n$ and
$r$ be positive integers such that  
$k$ contains a primitive $r$-th root of unity, $p\nmid n$, and 
$r/\!\gcd(n,r)$ divides $p^e-1$.
If both {\em(i)} and {\em(ii)} below are satisfied, then any
Galois extension $E/k(t)$ having ramification as in {\em ({$\dagger$})}
must have the form $E=k(v,w)$ where $v^{p^e}-v=\gamma w^n$ and 
$\gamma\in k^*$ and $t=w^r$.
\begin{enumerate}
\item $n$ is the least nonnegative integer congruent modulo~$r$
 to any number of the form $np^i$ with $i\ge 0$.
\item Either $k=\F_{p^e}$ or $p^e$ is the least power of $p$ which is
 congruent to $1$ modulo $r/\!\gcd(n,r)$.
\end{enumerate}
Conversely, if either {\em(i)} or {\em(ii)} is not satisfied, 
then there exist Galois extensions $E/k(t)$ having ramification
as in {\em ($\dagger$)} which do not have the above form.
\end{thm}
Condition (i) seems unexpected in this context.  It is especially
surprising that this condition is true for $n=(p^e+1)/4$ and $r=(p^e-1)/2$
(assuming $p^e\equiv 3\pmod{4}$ and $p^e>3$); in this case the subgroup of
$(\bZ/r\bZ)^*$ generated by $p$ has order $e$, but all $e$ elements of
the coset of $n$ have least nonnegative residue lying in the top half
of the interval $[0,r]$.

In fact, we do rather more than classify the two special types
of polynomials described above.  We determine all polynomials of
a general class which contains the polynomials of the two special 
types.  In particular, we prove the following result (see
Theorem~\ref{genpol2} for a refined version).
\begin{thm}
\label{genpol}
Let $k$ be an algebraically closed field of characteristic $p>0$,
let $d=(q^2-q)/2$ for some power $q=p^e$, and let 
$f(X)\in k[X]\setminus k[X^p]$ have degree $d$.
Then the following are equivalent:
\begin{enumerate}
\item $G:=\Gal(f(X)-u,k(u))$ has a transitive normal subgroup isomorphic to\/
$\PSL(q)$, and the Galois closure of the extension $k(x)/k(f(x))$
does not have genus $(q^2-q)/2$.
\item There exist linear polynomials $\ell_1,\ell_2\in k[X]$
   such that the composition $\ell_1\circ f\circ\ell_2$ is one
   of the following polynomials or one of the exceptions in Table~B
   (which follows Theorem~\ref{genpol2}):
\[X(X^m+1)^{(q+1)/(2m)}
 \left(\frac{(X^m+1)^{(q-1)/2}-1}{X^m}\right)^{(q+1)/m}\]
with $q$ odd and $m$ a divisor of $(q+1)/2$; or
\[X^{-q}\biggl(\sum_{i=0}^{e-1}X^{m2^i}\biggr)^{(q+1)/m}\]
with $q$ even and $m$ a divisor of $q+1$ with $m\neq q+1$.
\end{enumerate}
In these examples, $G\cong\PSL(q)$ if $m$ even, and
$G\cong\PGL(q)$ if $m$ odd.
\end{thm}
In the examples listed in (ii) (ignoring those in Table~B),
the cover $f\colon\Line\to\Line$ is only ramified over $\infty$ and $0$,
and any inertia group at a point over $0$ (in the Galois closure
cover) is cyclic of order $(q+1)/m$.  There is a point
over $\infty$ (in the Galois closure cover) whose inertia group
is the group of upper-triangular matrices in $G$, and whose
higher ramification groups (in the lower numbering) satisfy
$I_1=I_2=\dots=I_n\ne I_{n+1}=1$, where $n=m/\gcd(m,2)$.
The Galois closure of $k(x)/k(f(x))$ is $k(v,w)$ where
$v^q-v=w^n$.
We also prove some results in case the Galois closure
has genus $(q^2-q)/2$; we complete the analysis of this case
in the papers \cite{GRZ,GRZ2}.

This paper is organized as follows.  In the next section
we determine all group theoretic possibilities which could
correspond to a polynomial as in the previous theorem.  In Section~3 we
examine when the group theoretic data determines the Galois
closure $E$, and in particular prove Theorem~\ref{unq}.
In Section~4 we use knowledge of $E$ to classify polynomials satisfying
the properties discussed above.  For convenience, we collect various
elementary group theoretic facts in an appendix.

The second author thanks Hendrik Lenstra and Henning Stichtenoth
for valuable conversations.

\textbf{Notation:} Throughout this paper, $k$ is a field of
characteristic $p\ge 0$, and $\bar{k}$ is an algebraic closure
of $k$.  Also $q=p^e$ and $d=q(q-1)/2$. The letters $X$ and $Y$ denote
indeterminates, and (in situations where $k$ is present) the letters
$t,u,v,w,x,z$ denote elements of an extension of $k$ which are transcendental
over $k$.

\section{Group theory}

In this section we determine the possibilities for ramification in
the extension $k(x)/k(f(x))$, where $f(X)\in k[X]$ is a polynomial of
degree $q(q-1)/2$ whose arithmetic monodromy group $A$ has a transitive
normal subgroup isomorphic to $\PSL(q)$ (with $q$ being a power of
$\charp(k)$).  We denote the Galois closure of the extension 
$\bar{k}(x)/\bar{k}(f(x))$ by $E$, and write $g_E$ for its genus.
Then the geometric monodromy group of $f$ is
$G=\Gal(E/\bar{k}(f(x)))$, and we let
$G_1=\Gal(E/\bar{k}(x))$ denote
a one-point stabilizer of the permutation group~$G$.  Recall that, 
if $G$ is either $\PSL(q)$ or $\PGL(q)$, then a
Borel subgroup of $G$ is any subgroup
conjugate to the subgroup of upper-triangular matrices.
We often use without comment
the various elementary group theoretic facts collected in the Appendix.

\begin{thm} 
\label{groups}
Let $k$ be a field of characteristic $p>0$, and let $f(X)\in k[X]\setminus
k[X^p]$ have degree $d=q(q-1)/2$ where $q=p^e$.
If $A$ has a transitive normal subgroup $L$ isomorphic to\/ $\PSL(q)$, 
then all of the following hold unless $q,G,G_1$ are listed in Table~A\emph{:}
\begin{enumerate}
\item either $G=\PGL(q)$ or both $G=\PSL(q)$ and 
   $q\equiv 3\pmod{4}$;
\item $G_1\cap L$ is a dihedral group of order $2(q+1)/o$,
  where $o=\gcd(2,q-1)$; also $q\ge 4$;
\item the inertia group of a place of $E$ lying over
   the infinite place of $\bar{k}(f(x))$ is a Borel subgroup $I$ of~$G$; the
   higher ramification groups of this place satisfy
   $V:=I_1=I_2=\dots=I_n\gneqq I_{n+1}=1$.
\item $E/\bar{k}(f(x))$ has at most two finite branch points
  (i.e., ramified finite places of $\bar{k}(f(x))$); the
   possibilities are:
\begin{itemize}
\item One finite branch point, whose inertia group is cyclic of
order $|G\col L|(q+1)/(on)$ where $n\mid \left((q+1)/\!\gcd(4,q+1)\right)$ and $n<q+1$,
and $g_E= (q-1)(n-1)/2$.
\item No finite branch points, where $q\equiv 0\pmod{4}$ and $n=q+1$
    and $g_E=(q^2-q)/2$.
\item One finite branch point, with inertia group of order two and
    second ramification group trivial, where $q\equiv 0\pmod{4}$ and $n=1$
    and $g_E=(q^2-q)/2$.
\item Two finite branch points, both with inertia groups of order two,
    of which precisely one is contained in~$L$; here $q\equiv 1\pmod{2}$
    and $n=1$, and also $G=\PGL(q)$ and $g_E=(q^2-q)/2$.
\end{itemize}
\end{enumerate}
\end{thm} 

\begin{center}
$\begin{array}{ccc}
q&G&G_1 \\
\hline 
4&\PGammaL(q)&C_5\rtimes C_4\\
11&\PSL(q)&A_4\\
11&\PGL(q)&S_4\\
23&\PSL(q)&S_4\\
59&\PSL(q)&A_5
\end{array}
$
\end{center}
\begin{center}
Table~A
\end{center}
 
In this paper we will determine all polynomials having either the first or
second ramification possibilities (as well as all polynomials corresponding
to the situations in Table~A).  We will determine the polynomials in
the third and fourth possibilities in the papers \cite{GRZ}
and \cite{GRZ2}, respectively, and in the latter paper we will also
determine all polynomials of other degrees whose arithmetic monodromy
group has a transitive normal subgroup isomorphic to $\PSL(q)$.
We will prove in this paper that the
fourth possibility does not yield any exceptional polynomials, or any
indecomposable polynomials that decompose over an extension field.
However, it turns out that the third possibility does yield exceptional
polynomials.  At the end of this section we include a result giving
more details about this possibility, which we will need in \cite{GRZ}.

The remainder of this section is devoted to the proof of Theorem~\ref{groups}.
First, since $L\le A\le S_d$, we must have
$|\PSL(q)|\le|S_d|$, so $q\ge 4$.  Next we show that $L\le G$.
If $L\not\le G$ then (since $G\unlhd A$) the group $L\cap G$ is
a proper normal subgroup of the simple group $L$, and hence is trivial.
Since $L$ and $G$ normalize one another and intersect trivially, they
must commute.
Recall that the centralizer in $S_d$ of any transitive subgroup
has order at most $d$ \cite[Thm.~4.2A]{DM}.  But we have shown that
the centralizer of $G$ has order at least $|L|>d$, contradiction.
Thus $L\le G$.

We now describe the group-theoretic constraints implied by the
hypotheses of Theorem~\ref{groups}.  The transitive subgroup
$G\le S_d$ satisfies $\PSL(q)\cong L\unlhd G$, and $G_1$ is a
point-stabilizer of $G$.  We can identify the permutation representation
of $G$ with the action of $G$ on the set of left cosets of $G_1$ in $G$.
Thus, our hypothesis of transitivity of $L$ says that $LG_1=G$.

We use valuation theory to describe the further group-theoretic
conditions.  We identify the places of $\bar{k}(f(x))$ with 
$\bar{k}\cup\{\infty\}$, and we say that an element of $\bar{k}\cup\{\infty\}$
is a branch point of $f$ if the corresponding place is ramified in
$E/\bar{k}(f(x))$.
For a place $Q$ of $\bar{k}(f(x))$, let $P$ be a place of $E$ lying
over $Q$, and denote the ramification groups of $P/Q$ by
$I_0(Q)$, $I_1(Q)$, $\ldots$.  Different choices of $P$ yield
conjugate ramification groups; this ambiguity is irrelevant in what follows.
We use the following standard properties of these groups:
\begin{itemize}
\item  each $I_i$ is normal in $I_0$;
\item $I_1$ is the (unique) Sylow $p$-subgroup of $I_0$;
\item $I_0/I_1$ is cyclic;
\item $I_i=1$ for all sufficiently large $i$.
\end{itemize}
These properties imply that $I_1$ is the semidirect product of $I_0$
by a cyclic group of order $|I_0\col I_1|$.
Moreover, since $\infty$ is totally ramified in
$\bar{k}(x)/\bar{k}(f(x))$, we have $I_0(\infty)G_1=G$.
We write $I=I_0(\infty)$ and $V=I_1(\infty)$.  

By combining the Riemann-Hurwitz genus formula with
Hilbert's different formula~\cite[Prop.~IV.4]{Se},
we can express the genus $g_E$ in terms of the
ramification groups:
\[
2g_E-2 = -2d + \sum_Q |G\col I_0(Q)|\sum_{i\ge 0} (|I_i(Q)|-1).
\]
By combining this formula with the analogous one for the
extension $E/\bar{k}(x)$, we obtain a formula for the
(non-Galois) extension $\bar{k}(x)/\bar{k}(f(x))$.  We need some
notation to state this formula.
For a subgroup $H$ of $G$, let $\orb(H)$ be the number of orbits of $H$,
and let $\text{orb}'(H)$ be the number of orbits with length coprime to $p$.
Define the \emph{index} of $Q$ to be
\[
\ind Q=
\sum_{i\geq 0}\frac{d-\orb(I_i(Q))}{|I_0(Q)\col I_i(Q)|}.\]
Then our `Riemann-Hurwitz' formula for $\bar{k}(x)/\bar{k}(f(x))$ is
\[
2d-2=\sum_Q\ind Q.\]
We also use the inequality
\[\ind Q\geq d-\orb'(I_0(Q)),\]
which comes from the relation between the different exponent and 
ramification index~\cite[Prop.~III.13]{Se}.

Our first lemma shows that $G$ is contained in the automorphism
group $\PGammaL(q)$ of $L$ (see the appendix for information
about this group); we simultaneously prove that
$G_1\cap L$ usually has the desired shape.  

\begin{lemma} 
\label{zero}
$G$ is contained in\/ $\PGammaL(q)$.
Also, except for the final four cases listed in Table~A, the group
$L_1:=G_1\cap L$ is dihedral of order $2(q+1)/o$.
\end{lemma}

\begin{proof} 
As noted above, we have $q\ge 4$.
Since $|L\col L_1|=|G\col G_1|=d$,
the order of $L_1$ is  $2(q+1)/o$.
This numerical information severely limits the possibilities for
$L_1$; from Dickson's classification of subgroups of $L$
(Theorem~\ref{dickson}),
it follows that $L_1$ 
has the desired shape except possibly if $L_1$ is one of
the following groups:
$A_4$ if $q=11$; $S_4$ if $q=23$; $A_5$ if $q=59$.
In any case, $L_1=N_L(L_1)$ (since $L$ is simple
and $L_1$ is maximal unless $q=4,7,9,$ or 11, and these 
special cases are easily handled).  Hence $G_1=N_G(L_1)$.
Since $C_G(L)$ is a normal subgroup of $G$ contained
in $G_1$, it must be trivial, so indeed $G$ embeds in
the automorphism group $\PGammaL(q)$ of $L$.
Finally, the three exceptional possibilities for $L_1$ occur
with $q$ prime, so $G$ is either $L$ or $\PGL(q)$, and 
when $L_1$ is $S_4$ or $A_5$ we cannot have $G=\PGL(q)$ 
since there is no $G_1\leq G$ with $|G_1\col L_1|=2$.
\end{proof}

Henceforth assume $q,G,G_1$ do not satisfy any of the final
four possibilities in Table~A, and write $L_1:=G_1\cap L$.
Then Lemma~\ref{unique_action} says that our action of $G$ (on the
set of left cosets of $G_1$ in $G$)
is uniquely determined (up to equivalence) by $G$ and $q$, and does
not depend on the specific choice of $G_1$.

\begin{lemma}  
\label{borel}
$I$ is a Borel subgroup of $L$ or of\/ $\PGL(q)$,
unless $q,G,G_1$ are listed in Table~A.
\end{lemma}

\begin{proof}
We first show that $V$ is a Sylow $p$-subgroup of $L$ if $q>4$.
If $e\le 2$ then $\PGammaL(q)/L$ has order coprime to $p$, 
so $V$ is contained in a Sylow $p$-subgroup of $L$; but this Sylow
subgroup has order $q$, which implies the claim since $q\mid |V|$.
Now assume $e\ge 3$ and $q\ne 64$.  Then Zsigmondy's theorem implies
$p^e-1$ has a primitive prime divisor $s$.  Here $s$ is coprime
to $2e$, and hence to $|\PGammaL(q)\col L|$, so any $\mu\in I$ of order
$s$ must lie in $L$.  Now, $\mu$ acts on $V\cap L$ by conjugation,
with the identity element as a fixed point, and every other orbit having
size $s$ (since the centralizer in $L$ of an order-$p$ element of $L$ is a
Sylow $p$-subgroup of $L$, and hence does not contain $\mu$).  Thus
$|V\cap L|\equiv 1\pmod{s}$; since also $|V\cap L|$ divides $q$, 
primitivity of $s$ implies that $|V\cap L|$ is either $1$ or $q$.
We cannot have $|V\cap L|=1$, since in this case $V$ would embed in
$\PGammaL(q)/L$, but the order $eo$ of this group is not divisible by
$qo/2$.  Hence $V$ is a $p$-group containing a Sylow $p$-subgroup of
$L$; after replacing $I$ by a suitable conjugate, we may assume that $V\cap L$
consists of upper-triangular matrices, and that $V$ is generated by $V\cap L$
and a group of field automorphisms.  The conjugated
$\mu$ still lies in $I\cap L$ (since $L$ is normal in $\PGammaL(q)$), and
must be upper-triangular (since it normalizes $V\cap L$), so it does not
commute with any nontrivial field automorphism (by primitivity).  Thus,
as above, $|V|\equiv 1\pmod{s}$, so primitivity implies $|V|$ is a power
of $q$.  But $|\PGammaL(q)\col L|=eo$ is less than $q$, so $V$ is a Sylow
$p$-subgroup of $L$.  This conclusion is also true for $q=64$, as can be
shown by a similar argument using an order-$21$ element of $I$ (which must
lie in $L$).

Now, for $q>4$, we may conjugate to assume $V$ consists of upper triangular
matrices.  Then the normalizer $N$ of $V$ in $\PGammaL(q)$ is the semidirect
product of the group of all upper triangular matrices with the group of
field automorphisms.  Consider an element $\nu\in N$ whose
image in $\PGammaL(q)/\PGL(q)$ has order $r$; then the order of $\nu$ divides
$rp(q^{1/r}-1)$.  Since $I$ is generated by $V$ and an element $\nu$ of order
divisible by $(q-1)/o$, and $I\le N$, for this $\nu$ we have either
$r=1$ or $q=9$.  If $r=1$ then indeed $I$ is a Borel subgroup of either
$L$ or $\PGL(q)$.  If $q=9$ and $r>1$ then the above divisibility relations
imply $|I|=36$, so the conditions $IG_1=G$ and $|G\col G_1|=36$ imply $I\cap G_1=1$;
but one can check that there are no such $G$ and $G_1$.

Finally, for $q=4$, we have $\PGL(q)\cong A_5$ and $\PGammaL(q)\cong
S_5$, and one easily checks that, if $I$ is not a Borel subgroup of $L$,
we must have $G=\PGammaL(q)$ and $I\cong C_6$, with $G_1$ being the normalizer
in $S_5$ of a dihedral group of order 10.
\end{proof}

Henceforth we assume $q,G,G_1$ are not among the triples listed
in Table~A.  An immediate consequence of this lemma is the structure of the
higher ramification groups over the infinite place, as described
in (iii) of Theorem~\ref{groups}: recall that the chain
$I\geq V\geq I_2(\infty)\geq I_3(\infty)\geq\dots$ is such that
every $I_i(\infty)$ is normal in $I$, and $I_i(\infty)=1$ for some $i$.
In our case $V (=I_1(\infty))$ is a minimal normal subgroup of $I$, so 
every nontrivial $I_i(\infty)$ equals $V$; hence there is some $n\geq 1$ 
for which $V=I_1(\infty)=I_2(\infty)=\dots=I_n(\infty)$ but 
$I_{n+1}(\infty)=1$.

We complete the proof of Theorem~\ref{groups} by analyzing the possibilities
for the various other chains of ramification groups in light of
the Riemann-Hurwitz formula.  We do this in four cases.

\subsection{The case $q\equiv 3\pmod{4}$ and $I\leq L$}

First we compute the number of fixed points for various types of
elements of $\PGammaL(q)$ in the action under consideration.  We will
use these computations to determine the number and length of the
orbits of cyclic subgroups of $\PGammaL(q)$, which control the
indices of finite branch points for the
extension $\bar{k}(x)/\bar{k}(f(x))$.
\begin{lemma} 
\label{fps}
Let $\nu\in\PGammaL(q)$ have order $r>1$.  
For $\nu\in L$: if $r=2$ then
$\nu$ has $(q+3)/2$ fixed points; if $r$ divides $p(q-1)/2$ then
$\nu$ has no fixed points; if $r>2$ and $r$ divides $(q+1)/2$ then
$\nu$ has one fixed point.  If $\nu$
is a field automorphism then it has $(q_0^2-q_0)/2$ fixed points,
where $q_0=q^{1/r}$.  An involution in\/ $\PGL(q)\setminus L$
has $(q-1)/2$ fixed points.  
\end{lemma}
\begin{proof}
Since all $(q^2-q)/2$ involutions of $L$ are conjugate, and the
centralizer in $L$ of any such involution is dihedral of order $q+1$,
we can identify the action of $L$ on cosets of $L_1$ with the conjugation
action of $L$ on the involutions of $L$.  Then the uniqueness of the
action of $G$ on cosets of $G_1$ implies that that action is equivalent to the
conjugation action of $G$ on the involutions of $L$, so we examine
the latter action.  Thus, the number of fixed points of an element
$\nu\in G$ equals the number of involutions of $L$ which commute
with $\nu$.

If $\nu\in L$ is an involution, then its centralizer in $L$ is
dihedral of order $q+1$, and so contains $(q+3)/2$ involutions.
Conversely, if the order of $\nu\in L$ does not divide $q+1$,
then $\nu$ cannot lie in the centralizer of any involution of $L$.
If $\nu\in L$ has order $r>2$, where $r$ divides $(q+1)/2$, then
the normalizer of $\langle\nu\rangle$ in $L$ is dihedral of
order $q+1$, so the centralizer of $\nu$ in $L$ contains a unique
involution.
The centralizer in $L$ of an involution of $\PGL(q)\setminus L$
is dihedral of order $q-1$, and so contains $(q-1)/2$ involutions.
Finally, suppose $\nu\in\PGammaL(q)$ is a field automorphism of
order $r$, and let $q_0=q^{1/r}$; then the centralizer of $\nu$
in $L$ is $\PSL(q_0)$, which contains $(q_0^2-q_0)/2$ involutions.
\end{proof}

We compute $\ind\infty=(d-1)+n(q-1)$.  For any $Q\in \bar{k}$,
if $I_0(Q)\neq 1$ then $I_0(Q)$ contains an element $\nu$ of prime
order $r$;
by Lang's theorem, this element is either in $\PGL(q)$ or is
conjugate to a field automorphism.  If $\nu$ is an involution
of $L$, then $\ind Q\geq d-\text{orb}'(I_0(Q))\geq
d-\text{orb}'(\langle\nu\rangle)=(d-(q+3)/2)/2$.  If $\nu\in L$ is
not an involution, then likewise $\ind Q\geq d-(d/r)$ if $d\mid
(q-1)$, and $\ind Q\geq d - 1-(d-1)/r$ if $r\mid (q+1)$, and 
$\ind Q\geq
d$ if $r=p$.  If $\nu$ is an involution in $\PGL(q)\setminus L$
then $\ind Q\geq (d-(q-1)/2)/2$.  Finally, if $\nu$ is a field
automorphism then $\ind Q\geq d-s-(d-s)/r$, where $s=(q_0^2-q_0)/2$
and $q_0^r=q$;
moreover, if the field automorphism has order $p$ then $\ind Q\geq d-s$.
In any case, we conclude that 
$\ind Q\geq (d-(q+3)/2)/2.$   

If there are no finite branch points, then $(d-1)+n(q-1)=2d-2$,
so $q-1\mid d-1$ which is absurd; hence there is at least one finite
branch point.
If there are at least two finite branch points, then 
\[(d-1)+n(q-1)+\sum_{Q}\ind Q\geq
(d-1)+(q-1)+(d-(q+3)/2)>2d-2,\] 
contrary to the Riemann-Hurwitz formula.  Hence there is exactly one
finite branch point~$Q$.
Similarly, the Riemann-Hurwitz formula would be violated if $I_0(Q)$
contained an element of order $p$; thus, $p\nmid |I_0(Q)|$.

Now we use another geometric fact---equivalently, we add another
condition to the list of group theoretic restrictions.  Consider the
field extension
$E^L/E^G$: it is unramified over the infinite place
of $E^G$ (since $L$ contains all conjugates of $I$), so it is 
ramified over at most the 
single place $Q$ of $E^G$, and the ramification over that place is
tame.  Since any nontrivial tame extension of $E^G=\bar{k}(f(x))$ ramifies
over at least two places, we must have $L=G$.

Next we consider the extension $E/E^{L_1}$, which is Galois
with group dihedral of order $q+1$; here $E^{L_1}=\bar{k}(x)$ has
genus zero.
This extension is unramified over the infinite place
of $E^{L_1}$ (since $\gcd(|L_1|,|I|)=1$), so all ramification occurs over
places of $\bar{k}(x)$ lying over $Q$.   
Let $H$ be the cyclic subgroup of $L_1$ of order $(q+1)/2$;
then $E^H/E^{L_1}$ has degree 2 and ramifies only over places
of $\bar{k}(x)$ lying over $Q$.  But $E^H/E^{L_1}$ is a nontrivial
extension of $\bar{k}(x)$, so it is ramified, whence 
the ramification index of $Q$ (in $E/E^G$) is even.

Once again considering the extension $E/E^G$, we have shown
that there is exactly one finite branch point $Q$, which is tamely
ramified of even index.  In particular, the inertia group $I_0(Q)$ is
cyclic of order $2s$, where $p\nmid s$.  Since $I_0(Q)\leq L$, by
Theorem~\ref{dickson} we have $s\mid (q+1)/4$.  By Lemma~\ref{fps},
$I_0(Q)$ has $(d-(q+3)/2)/(2s)$ orbits of length $2s$; if $s=1$ then
the other points are fixed, while if $s>1$ then only one point is
fixed and $I_0(Q)$ has $(q+1)/(2s)$ orbits of length $s$.
In any case, $\ind Q=d-1-(q^2-1)/(4s)$.
{}From $2d-2=(d-1)+n(q-1)+\ind Q$ we conclude that $n=(q+1)/(4s)$.

Now we compute the genus of $E$, using the Riemann-Hurwitz
formula for the extension $E/E^G$: namely, $2g_E-2$ equals
\begin{align*}
-2&|G|+|G\col I|((|I|-1)+n(|V|-1))+|G\col I_0(Q)|(2s-1) \\
&=-(q^3-q)+(q+1)\left(
\frac{q^2-q}{2}-1+n(q-1)\right)+\frac{q^3-q}{4s}(2s-1);
\end{align*}
simplifying, we find that $g_E=(n-1)(q-1)/2$. 
 
This concludes the proof of Theorem~\ref{groups} in this case.

\subsection{The case $q\equiv 3\pmod{4}$ and $I\not\leq L$}

Here $I$ is a Borel subgroup of $\PGL(q)$, so $G$ contains
$\PGL(q)$.  As above, we compute $\ind\infty=d-1+n(q-1)/2$, and each
finite branch point $Q$ satisfies $\ind Q\geq
(d-(q+3)/2)/2$; hence there are at most two finite branch
points.  

If there are two finite branch points $Q_1$ and $Q_2$ then 
\[d-1+n(q-1)/2+\sum_{Q}\ind Q\geq d-1+(q-1)/2+(d-(q+3)/2)=2d-3,\]
so we must have $n=1$, and the indices $\ind Q_1$ and $\ind Q_2$ must
be (in some order) $(d-(q+3)/2)/2$ and
$(d-(q-1)/2)/2$.  The inequalities for $\ind Q$ in the previous
subsection imply that $I_0(Q_1)$ and $I_0(Q_2)$ are 2-groups containing no
involutions outside $\PGL(q)$.
Since these two inertia groups have order coprime to $p$, they 
are cyclic, so they must have order 2 (otherwise $\ind{Q_j}$
would be too large).  
Assuming without loss that $\ind Q_1<\ind Q_2$, it
follows that $I_0(Q_1)$ is generated by an
involution in $L$ and $I_0(Q_2)$ is generated by an involution in
$\PGL(q)\setminus L$.  The subgroup of $G$ generated by the
conjugates of $I_0(Q_1)$ and $I_0(Q_2)$ is $\PGL(q)$, so
$E^{\PGL(q)}/E^G$ is an unramified extension of
$\bar{k}(f(x))$, and thus is trivial: $G=\PGL(q)$.  Now we compute the
genus of $E$:
\[2g_{E}-2
 = -2(q^3-q)+(q+1)(q^2-q-1+q-1)+(q^3-q),\]
so $g_{E}=(q^2-q)/2$. 

Henceforth assume there is at most one finite branch point $Q$.
The same argument as in the previous subsection shows there must be
at least one such $Q$, and it must be tamely ramified.
Consider the extension 
$E^{\PGL(q)}/E^G$.  This extension is unramified over the
infinite place of $\bar{k}(f(x))$ (since $\PGL(q)$ contains all conjugates
of $I$), so it can only
be ramified over $Q$; hence it is a tamely ramified extension
of $\bar{k}(f(x))$ which is ramified over less than two points, so
it is trivial: $G=\PGL(q)$.  Thus $E^L/E^G$ has degree 2,
so it is a tame extension of $E^G=\bar{k}(f(x))$, whence it must be
ramified over both infinity and $Q$; hence $I_0(Q)$ is cyclic of
order $2s$ and is not contained in $L$.
It follows that $2s$ divides either $q-1$ or $q+1$.  In the former
case, a Riemann-Hurwitz computation shows that $n=1$ and
$s=(q-1)/2$, but then the genus of $E$ would be $-q$ which is
absurd.  In the other case, we find that $s$ must be even, 
$n=(q+1)/(2s)$, and the genus of $E$ is $(q-1)(n-1)/2$.
This concludes the proof of Theorem~\ref{groups} in case $q\equiv 3$~(mod~4).

\subsection{The case $q\equiv 1$~(mod~4)}

First we show that $BL_1\ne L$ for every Borel subgroup $B$ of $L$; this
implies $B\ne I$, so $I$ is a Borel subgroup of $\PGL(q)$.
For, if $BL_1=L$, then every Borel subgroup of $L$ has the form $B^{\nu}$
with $\nu\in L_1$, so $B^{\nu}\cap L_1=(B\cap L_1)^\nu$ is trivial.
But this is impossible, since any involution in $L_1$ is contained in some
Borel subgroup of $L$.

Thus, $I$ is a Borel subgroup of $\PGL(q)$,
so $G\geq \PGL(q)$.  The remainder of the proof in this case is
similar to the proof in the previous subsection, so we only give
the fixed point computation.
\begin{lemma}
Let $\nu\in\PGammaL(q)$ have order $r$.  If $\nu$ is
a field automorphism then $\nu$ has no fixed points if $r$ is
even, and $\nu$ has $(q_0^2-q_0)/2$ fixed points if $r$ is odd; 
here $q_0=q^{1/r}$.  If 
$\nu\in \PGL(q)$ and $r>2$, then $\nu$ has one fixed point
if $r\mid(q+1)$ and $\nu$ has no fixed points if $r\mid p(q-1)$.
An involution in $L$ has $(q-1)/2$ fixed points; an involution
in\/ $\PGL(q)\setminus L$ has $(q+3)/2$ fixed points.  
\end{lemma}
The proof of this lemma is similar to the proof of Lemma~\ref{fps},
except that in this case we identify the action of $G$ on cosets of $G_1$
with the
conjugation action of $G$ on the involutions in $\PGL(q)\setminus L$.

\subsection{The case $p=2$}

Here $\PGL(q)\cong L$, so we identify these groups.  In particular,
$I$ is a Borel subgroup of $L$.
We begin with a fixed point calculation.  
\begin{lemma}
An involution of $L$ has $q/2$ fixed
points; an element of $L$ of order $r>1$ has no fixed points if
$r\mid q-1$, and one fixed point if $r\mid q+1$.  A
field automorphism of order $r$ has no fixed points if $r=2$,
and has $(q_0^2-q_0)/2$ fixed points if $r$ is an odd prime; here 
$q_0=q^{1/r}$.
\end{lemma}
\begin{proof}
Since there are $(q^2-q)/2$ dihedral subgroups of $L$ of order
$2(q+1)$, and they are all conjugate, the number of fixed points of
an element $\nu\in L$ (acting on cosets of $L_1$ in $L$) equals the number of 
dihedral subgroups of $L$ of order $2(q+1)$ which contain $\nu$.
This number is certainly zero if the order $r>1$ of $\nu$ does not
divide $2(q+1)$, which happens if $r\mid q-1$.
Next, a dihedral group of order $2(q+1)$ normalizes each of its
cyclic subgroups of order dividing $q+1$; since the normalizer in $L$
of a nontrivial cyclic group of order dividing $q+1$ is dihedral
of order $2(q+1)$, it follows that a nonidentity element of $L$ of
order dividing $q+1$ has exactly one fixed point.

Next consider a field automorphism $\sigma$ of order 2.  Let $\bar{G}$
be the group generated by $L$ and $\sigma$, so $|\bar{G}\col L|=2$.
Suppose that $\sigma$ has a fixed point; then $\sigma$ lies in some
point-stabilizer $\bar{G}_1\leq\bar{G}$.  Here $\bar{L}_1:=\bar{G}_1\cap L$
is dihedral of order $2(q+1)$, and $\bar{G}_1$ is the normalizer of
$\bar{L}_1$ in $\bar{G}$. In particular, $\sigma$ normalizes the 
cyclic subgroup of $\bar{L}_1$ of order $q+1$; let $\mu$ be a generator
of this subgroup.  Let $\hat\mu\in\GL(q)$ be a preimage of $\mu$ under
the natural map $\GL(q)\to L$.  Then $\hat\mu$ has 
distinct eigenvalues, and the eigenvalues of $\sigma \hat\mu\sigma^{-1}$
are the $\sqrt{q}$-th powers of the eigenvalues of $\hat\mu$, so 
$\sigma \hat\mu\sigma^{-1}$ is conjugate to $\hat\mu^{\sqrt{q}}$.  Thus, 
$\sigma \mu\sigma^{-1}$ is conjugate to $\mu^{\sqrt{q}}$ in $L$.
Since the normalizer of $\langle \mu\rangle$ in $L$ is $\bar{L}_1$, the
only powers of $\mu$ which are conjugate (in $L$) to $\mu^{\sqrt{q}}$
are $\mu^{\pm\sqrt{q}}$.  It follows that $\sigma \mu\sigma^{-1}$ is one
of $\mu^{\pm\sqrt{q}}$.  A Sylow $2$-subgroup of $\bar{G}_1$ which contains
$\sigma$ also contains an involution $\tau\in\bar{L}_1$, and so must be
dihedral of order 4; thus $(\sigma \tau)^2=1$.  But, since $\tau \mu\tau=\mu^{-1}$,
we have
$\mu=(\sigma \tau)^2 \mu(\sigma \tau)^{-2}=\mu^q=\mu^{-1}$,
which is a contradiction, so in fact $\sigma$ has no fixed
points.  

For the remaining cases, we use the following elementary result about fixed
points in a permutation group.  Let $\bar{G}_1$ be a subgroup of the
finite group $\bar{G}$, and consider the natural action of 
$\bar{G}$ on the set of cosets of $\bar{G}_1$ in $\bar{G}$.  For an element $\nu\in \bar{G}$,
let $\nu_1,\nu_2,\dots,\nu_s$ be representatives of the conjugacy
classes of $\bar{G}_1$ which are contained in the conjugacy class of
$\nu$ in $\bar{G}$.  Then the number of fixed points of $\nu$ is 
\[\sum_{i=1}^s |C_{\overline{G}}(\nu)\col C_{\overline{G}_1}(\nu_i)|.\]

We first apply this in case $\nu\in L$ is an involution; let
$\bar{G}=L$ and let $\bar{G}_1=L_1$ be dihedral of order $2(q+1)$. 
Then the conjugacy class of $\nu$ in $L$ consists of all involutions in
$L$, and $C_L(\nu)$ is a Sylow $2$-subgroup of $L$; all involutions in
a dihedral group of order $2(q+1)$ are conjugate, and each generates
its own centralizer.  By the above fixed point formula, $\nu$ fixes $q/2$
cosets of $L_1$ in $L$.

For the final case in the lemma, let $\sigma\in\PGammaL(q)$ be a
field  automorphism of order $r$, and
assume that $r$ is an odd prime.  Write $q=q_0^r$.  Let $\bar{G}$ be the group
generated by $L$ and $\sigma$, so $|\bar{G}\col L|=r$.  Let $H$ be a 
cyclic subgroup of $\PSL(q_0)$ of order $q_0+1$, and
let $\bar L_1$ (respectively, $\bar G_1$) be the normalizer of $H$ in 
$L$ (respectively, $\bar{G}$).  Then $\bar L_1$ is
dihedral of order $2(q+1)$, so
$\bar G_1$ is generated by $\bar L_1$ and $\sigma$ and has order
$2r(q+1)$.  The centralizer of $\sigma$ in $\bar G_1$ is generated by
$\sigma$ and $\PSL(q_0)\cap{\bar G_1}$; since the latter group is dihedral
of order $2(q_0+1)$, the order of the centralizer is $2r(q_0+1)$.
The centralizer of $\sigma$ in $\bar{G}$ is
generated by $\sigma$ and $\PSL(q_0)$, and has order $r(q_0^3-q_0)$. 
We will show that the conjugacy class of $\sigma$ in $\bar{G}$
contains a unique conjugacy class of $\bar{G}_1$; by the fixed point
formula, it follows that $\sigma$ has $(q_0^2-q_0)/2$ fixed points.

Since $\sigma$ normalizes $L$, any conjugate of $\sigma$ in $\bar{G}$ has
the form $\mu\sigma$ with $\mu\in L$; conversely, by Lang's
theorem, if $\mu\sigma$ has order $r$ (where $\mu\in L$) then
$\mu\sigma$ is conjugate to $\sigma$ in $\bar{G}$.
If $r\nmid q+1$ then $\sigma$ generates a Sylow subgroup of
$\bar{G}_1$, so it is conjugate in $\bar{G}_1$ to any element 
$\mu\sigma$ of order $r$ with $\mu\in\bar{L}_1$.  Henceforth
assume $r\mid q+1$.  Let $\bar{J}$ be a Sylow $r$-subgroup of $\bar{L}_1$ 
containing $\sigma$.
Then $J:=\bar{J}\cap L$ is a Sylow $r$-subgroup of $\bar{L}_1$, so it
is cyclic of order (say) $r^j$; let $\tau$ be a generator of $J$.
Since $\sigma$ normalizes $J$, we have $\sigma \tau\sigma^{-1}=\tau^i$
where $1\leq i<r^j$.  Then $\tau^{i^r}=\sigma^r \tau\sigma^{-r}=\tau$,
so $i^r\equiv 1\pmod{r^j}$, whence $i\equiv 1\pmod{r^{j-1}}$.
Since $q+1=q_0^r+1$ is divisible by $r$, also $r\mid q_0+1$, 
and thus $(q+1)/(q_0+1)$ is divisible by $r$ but not by $r^2$.
Hence $\sigma$ does not centralize $J$, so $i>1$.
Next, we compute $(\tau^a\sigma)^r=\tau^{a(1+i+\dots+i^{j-1})}$,
so $\tau^a\sigma$ has order $r$ if and only if $r^j$ divides
$a(1+i+\dots+i^{j-1})$, i.e., if and only if $r^{j-1}\mid a$.
Conversely, if $r^{j-1}\mid a$ then $a\equiv b(i-1)\pmod{r^j}$;
since $\tau^{-b}\sigma \tau^b=\tau^{b(i-1)}\sigma$, it follows
that $J$ contains a unique conjugacy class of subgroups of order $r$
which are not contained in $L$.
Thus $\bar{G}_1$ contains a unique conjugacy class of subgroups of
order $r$ which are not contained in $L$, and the proof is complete.
\end{proof}
 
This computation implies that the index of
any branch point $Q$ satisfies $\ind Q> 2d/3-q/2$.
Since $\ind\infty=d-1+n(q/2-1)\geq d-2+q/2$, it follows that there is at
most one finite branch point $Q$, and $\ind Q\leq d-(q/2)$.

Since $G/L$ is cyclic, there is a unique group $H$ between $L$ and
$G$ such that $|H\col L|$ is the highest power of 2 dividing $|G\col L|$.
Then $E^H/E^G$ is unramified over infinity, so it is a 
tame extension of $\bar{k}(f(x))$
having only one branch point, whence it is trivial.  Thus $H=G$, so
if $G\neq L$ then $|G\col L|$ is a power of 2.  Now assume $s:=|G\col L|$ is
a power of 2.  Since $G$ is generated by the inertia groups, which
are conjugates of $I$ and $I_0(Q)$, it follows that $I_0(Q)$ maps
onto $G/L$; let $\nu\in I_0(Q)$ map to a generator of $G/L$.
Replacing $\nu$ by an odd power of itself, we may assume $\mu$
has order a power of 2.  Since $\ind Q\leq d-(q/2)$, $\nu$ must
have at least $q/2$ fixed points; hence $\nu$ cannot be conjugate
to a field automorphism, so $\langle\nu\rangle$ intersects $L$
nontrivially.  Since $4\nmid |L_1|$, a subgroup of $L$ of order
divisible by 4 has no fixed points; thus $4\nmid |I_0(Q)\cap L|$.
Hence the Sylow $2$-subgroup $I_1(Q)$ of $I_0(Q)$ has order $2s$, so it
is generated by~$\nu$.  Since the involution in $I_1(Q)\cap L$ is
centralized by $I_0(Q)\cap L$, the latter group has order a
power of 2; hence $I_0(Q)=I_1(Q)$.
By Lemma~\ref{serre}, $I_1(Q)/I_2(Q)$ is elementary abelian.
Since $I_1(Q)$ is cyclic of order $2s$, the order of $I_2(Q)$ is
either $e$ or $2e$, so $|I_2(Q)|>1$.
Thus, for $0\leq i\leq 2$, the group $I_i(Q)$ contains an involution of $L$, so
$\orb(I_i(Q))\leq q^2/4$.  But then we have the contradiction
\begin{eqnarray*}
\ind Q&\geq& \sum_{i=0}^2 (d-\orb(I_i(Q)))/|I_0(Q)\col I_1(Q)|\\
     &\geq& (d-q^2/4)+(d-q^2/4)+(d-q^2/4)/2> d-q/2.
\end{eqnarray*}
Henceforth we assume $G=L$.
 
Suppose the ramification over $Q$ is wild.  If $I_0(Q)$ has a 
subgroup of order 4 then every orbit of $I_0(Q)$ has even length, 
so $\ind Q\geq d$ is too large.  Hence the Sylow $2$-subgroup $I_1(Q)$
of $I_0(Q)$ has order 2, so it is centralized by $I_0(Q)$; but
the centralizer (in~$L$) of an involution is the Sylow subgroup
 containing the involution, so $I_0(Q)$ has order 2.  A
Riemann-Hurwitz calculation yields $I_2(Q)=1$ and $n=1$.
In this case the genus of $E$ is $(q^2-q)/2$.

Now suppose the ramification over $Q$ is tame.  Then $I_0(Q)$ is
cyclic of order $s$, where $s$ divides either $q-1$ or $q+1$.
If $s>1$ and $s\mid q-1$ then $\ind Q=d-d/s$, which implies
that $n=1$ and $s=q-1$; but then the genus of $E$ would 
be $-q$, which is absurd.  Thus $s\mid q+1$, so 
$\ind Q=d-1-(d-1)/s$, whence $n=(q+1)/s$.   Here the genus 
of $E$ is $(q-1)(n-1)/2$.  This completes the proof of Theorem~\ref{groups}.

\subsection{Preliminary analysis in the case of two wildly ramified branch points}

We now prove the following result, which will be used in \cite{GRZ}:

\begin{lemma}
\label{grz}
Let $k$ be a field of characteristic $2$, and let $q=2^e$ with $e>1$.
Let $f\in k[X]\setminus k[X^2]$ have degree
$q(q-1)/2$, let $E$ be
the Galois closure of $k(x)/k(f(x))$, and let $\ell$ be the algebraic closure
of $k$ in $E$.  Suppose that $\Gal(E/\ell(f(x)))\cong\PGL(q)$ and
$k(x)/k(f(x))$ is wildly
ramified over at least two places of $k(f(x))$.  Then $E/\ell(f(x))$ has
precisely two ramified places, both of degree one, and the corresponding
inertia groups are
(up to conjugacy) the order-$2$ group generated by $\textmatrix{1&1\\0&1}$
and the group of upper-triangular matrices in\/ $\PGL(q)$.  Moreover, the second
ramification group over each ramified place is trivial.  The degree $[\ell\col k]$
divides $e$, and $f$ is indecomposable.  Here $f$ is exceptional if and only if
$e$ is odd and $[\ell\col k]=e$.
Finally, there is a curve $\Cone_0$ over $k$ such that $\ell.k(\Cone_0)\cong_{\ell} E$.
\end{lemma}

\begin{proof}
The facts about ramification follow from Theorem~\ref{groups}, using
the fact that $\PGL(q)$ contains a unique conjugacy class of involutions
(Lemma~\ref{involutions}).  Write $u:=f(x)$, and let $A:=\Gal(E/k(u))$ and
$G:=\Gal(E/\ell(u))$ and $G_1:=\Gal(E/\ell(x))$.  Theorem~\ref{groups} implies
$G_1$ is dihedral of order $2(q+1)$, so $G_1$ is maximal in $G$.
Pick $\nu\in C_A(G)$.  The subfield of $E$ fixed by
$G.\langle \nu\rangle$
is $\ell^{\nu}(u)$.  Since $\Gal(\ell(x)/\ell^{\nu}(x))\cong
\Gal(\ell(u)/\ell^{\nu}(u))\cong\Gal(\ell/\ell^{\nu})$, it follows that
$\Gal(E/\ell^{\nu}(x))$ contains an element
$\mu\nu$ with $\mu\in G$.  Since $\ell(x)/\ell^{\nu}(x)$ is Galois,
$\mu\nu$ normalizes $G_1$;
but since $\nu$ commutes with $G_1$, this means that $\mu$ normalizes $G_1$,
whence $\mu\in G_1$ (since $G$ is simple and $G_1$ is maximal in $G$).
Hence $\Gal(E/\ell^{\nu}(x))$
contains $\nu$, so $E^{\nu}$ contains
$\ell^{\nu}(x)$.  Since $\nu$ commutes with $G$, the group
$\langle \nu\rangle$ is a normal
subgroup of $G.\langle \nu\rangle$, so $E^{\nu}/\ell^{\nu}(u)$
is Galois.  But $E$ is the
splitting field of $f(X)-u$ over $\ell(u)$, so it is also the splitting field
of $f(X)-u$ over $\ell^{\nu}(u)$, whence it is the minimal Galois extension
of $\ell^{\nu}(u)$ which contains $\ell^{\nu}(x)$.
Thus $E=E^{\nu}$, so $\nu=1$.

We have shown that $C_A(G)=1$, so the action of $A$ on $G$ by conjugation is faithful, whence
$A\hookrightarrow \Aut(G)\cong\PGammaL(q)$.  In particular, $[\ell\col k]=|A\col G|$
divides $|\PGammaL(q)\col\PGL(q)|=e$.
Next, $f$ is indecomposable over $\ell$, since $G_1$ is a maximal subgroup
of $G$.  The existence of a curve $\Cone_0$ over $k$ with $\ell.k(\Cone_0)\cong_{\ell} E$
follows from the fact that $A$ is the semidirect product of $G$ with a group
$H$ of field automorphisms:
for $\ell^H=k$ and $|H|=|\Gal(\ell/k)|$, so $\ell.E^H=E$.

Finally, we address exceptionality.  By \cite[Lemma~6]{Coh} (or
\cite[Lemma~4.3]{GTZ}),
$f$ is exceptional if and only if every element of $A$ which generates $A/G$
has a unique fixed point.  By Theorem~\ref{groups}, $G_1$ is dihedral of order
$2(q+1)$; by Lemma~\ref{dihedrals}, $G_1$ is the normalizer of its unique subgroup of
order $3$, and $G$ has a unique conjugacy class of order-$3$ subgroups.
Thus we can identify the set of cosets of $G_1$ in $G$ with the set of order-$3$
subgroups of $G$, and the action of $G$ on cosets of $G_1$ corresponds to its action by
conjugation on these subgroups.  Moreover, the natural action of $A$ on the order-$3$
subgroups of $G$ induces the permutation action of $A$ under consideration.

Every coset in $A/G$ has the form $\sigma G$ where $\sigma\in\PGammaL(q)$ is a
field automorphism.  If $\sigma$ has order $e/e'$, then the centralizer of $\sigma$
in $G$ is $\PGL(2^{e'})$, which contains $2^{e'-1}(2^{e'}-1)$ subgroups of order $3$.
In particular, if $e'>1$ then $\sigma$ fixes more than one order-$3$ subgroup;
thus, if $f$ is exceptional then $A=\PGammaL(q)$.  Moreover, Theorem~\ref{fgsthm}
implies that $e$ must be odd if $f$ is exceptional.

Conversely, suppose $A=\PGammaL(q)$ and $e$ is odd.  To prove $f$ is exceptional,
we must show that every element of $A$ which generates $A/G$ has precisely one
fixed point.  By an easy counting argument \cite[Lemma~13.1]{FGS}, it suffices to
show that there is a generating coset of $A/G$ in which every element has at most
one fixed point.  Pick an element $\nu\in A$ which induces the Frobenius automorphism
on $A/G\cong\Gal(\F_q/\F_2)$.  Since $A/G$ has odd order,
$\nu^2$ also generates $A/G$, and moreover if $\nu^2$ has at most one fixed point
then so does $\nu$.  By replacing $\nu$ by $\nu^{2^s}$, where $s$ is sufficiently large
and $2^s\equiv 1\pmod{e}$, we may assume that $\nu$ has odd order; thus, if $\nu$ normalizes
an order-$3$ subgroup then it centralizes the subgroup.  By Lang's theorem on algebraic
groups, there exists $\tau\in\PGL(\bar{\F}_q)$ such that $\tau\nu\tau^{-1}$ is the Frobenius
automorphism $\sigma$ in $\PGammaL(\bar{\F}_q)$.  Thus $C_G(\nu)$ is isomorphic to a subgroup of
$C_{\PGL(\overline{\F}_q)}(\sigma)\cong\PGL(2)$, which contains a unique order-$3$ subgroup.
\end{proof}

The following consequence of Lemma~\ref{grz} describes the ramification in
$\Cone\to\Cone/B$, where $\Cone=\Cone_0\times_k \ell$ and $B$ is the group of upper-triangular
matrices in $L$.  Here $W=B\cap\PGL(2)$ and $T$ is the group of diagonal matrices in $B$.

\begin{cor}
If $\Cone$ is a curve over $\ell$ for which $\ell(\Cone)=E$, then the following hold:
\begin{enumerate}
\item $B$ acts as a group of $\ell$-automorphisms on $\Cone$;
\item the quotient curve $\Cone/B$ has genus zero;
\item the cover $\Cone \to \Cone/B$ has exactly three branch points;
\item the inertia groups over these branch points are $B$, $T$, and
      $W$ (up to conjugacy); and
\item all second ramification groups in the cover $\Cone \to \Cone/B$ are 
      trivial.
\end{enumerate}
\end{cor}

\begin{proof}
We know that $E/\ell(f(x))$ is Galois with group $L$, and is ramified over
precisely two places, both of which have degree one and have trivial second
ramification group, and moreover the inertia groups are $B$ and $W$.
Thus there is an action of $L$ as a group of $\ell$-automorphisms of $\Cone$
for which the cover $\Cone\to\Cone/L$ has this same ramification.
By Riemann-Hurwitz, the genus of $\Cone$ is $q(q-1)/2$.  Since the second
ramification groups in $\Cone\to\Cone/L$ are trivial, it follows that the
same is true in $\Cone\to\Cone/B$.

Let $P_1$ and $P_2$ be points of $\Cone$ whose inertia groups in
$\Cone\to\Cone/L$ are $B$ and $W$, and let $Q_1$ and $Q_2$ be the points
of $\Cone/L$ lying under $P_1$ and $P_2$.
Then each of the $q+1$ points of $\Cone$ lying over $Q_1$ has inertia group
conjugate to $B$.  Since $B$ is self-normalizing in $L$, and the intersection
of any two conjugates of $B$ is cyclic of order $q-1$, it follows that the
inertia group of $P_1$ in $\Cone\to\Cone/B$ is $B$, while the inertia groups
of the other $q$ points over $Q_1$ are cyclic of order $q-1$.  Since $B$ has
$q$ cyclic subgroups of order $q-1$, and they are all conjugate in $B$, it
follows that $Q_1$ lies under two branch points of $\Cone\to\Cone/L$, and
the corresponding inertia groups are $B$ and $T$.
Since all involutions in $L$ are conjugate, and each has normalizer of
order $q$, each of the $q^2-1$ involutions in $L$ occurs as the inertia
group of precisely $q/2$ points of $\Cone$ lying over $Q_2$.  Thus $Q_2$
lies under $q(q-1)/2$ points of $\Cone$ which ramify in $\Cone\to\Cone/B$,
and all these points lie over the same point of $\Cone/B$.  Hence
$\Cone\to\Cone/B$ has precisely three branch points, with inertia groups
$B$, $T$, and $W$; now Riemann-Hurwitz implies that $\Cone/B$ has genus
zero, which completes the proof.
\end{proof}

\section{Characterizing certain field extensions by ramification}

In this section we study the extensions $E/k(z)$ having
certain ramification; the specific choice of ramification 
configuration comes from our desired application to the classification of
degree $q(q-1)/2$ polynomials having monodromy group normalizing
$\PSL(q)$.  However, the results in this section
are of interest for their own sake, as they provide data for
the general problem of determining the ramification possibilities
for covers of curves having specified monodromy group.
The goal of this section is to prove Theorem~\ref{unq} of the introduction.

We begin with three easy facts about elementary abelian field
extensions; a convenient reference for
these is~\cite{GaSt}.  The first result describes the shape of these 
extension fields~\cite[Prop.~1.1]{GaSt}.

\begin{lemma} 
\label{one}
 Let $E/F$ be an elementary abelian extension of degree 
$p^e$, where\/ $\F_{p^e}\subseteq F$.  Then there exist elements 
$v\in E$ and $z\in F$ such that $E=F(v)$ and the minimal polynomial
of $v$ over $F$ is $T^{p^e}-T-z$.
\end{lemma}
The converse of this result also holds~\cite[Prop.~1.2]{GaSt}.

\begin{lemma}
\label{two}
Pick $z\in F$, and suppose that\/ $\F_{p^e}\subseteq F$ and that the
polynomial $T^{p^e}-T-z\in F[T]$ is irreducible over~$F$.  For any
root $v$ of this polynomial, the extension $F(v)/F$ is elementary abelian of
degree $p^e$.  The map $\sigma\mapsto\sigma(v)-v$ is an isomorphism
between $\Gal(F(v)/F)$ and the additive group of\/ $\F_{p^e}$.
The intermediate fields $F\subset F_0\subseteq F(v)$
with $[F_0\col F]=p$ are precisely the fields $F_0=F(v_\zeta)$, where for
$\zeta\in\F_{p^e}^*$ we put
\[v_\zeta:=(\zeta v)^{p^{e-1}}+(\zeta v)^{p^{e-2}}+\dots+(\zeta v)^p+
  (\zeta v).\]
The minimal polynomial for $v_\zeta$ over $F$ is $T^p-T-\zeta z$.
\end{lemma}

Note that the values $v_\zeta$ are precisely the images of the
various $\zeta v$ under the polynomial $T^{p^{e-1}}+\dots+T^p+T$,
which is the trace map from $\F_{p^e}$ to $\F_p$.
In order to apply this result, we need to know when the polynomial 
$T^{p^e}-T-z\in F[T]$ is irreducible; we now give criteria for
this~\cite[Lemma~1.3]{GaSt}.

\begin{lemma} 
\label{three}
Pick $z\in F$, and suppose\/ $\F_{p^e}\subseteq F$.  The 
following conditions are equivalent: \hfil\break
\item {\rm(a)} $T^{p^e}-T-z$ is irreducible over $F$.
\item {\rm(b)} For all $\zeta\in\F_{p^e}^*$, the polynomial $T^p-T-\zeta z$ is irreducible
  over~$F$.
\item {\rm(c)} For all $\zeta\in\F_{p^e}^*$, the polynomial $T^p-T-\zeta z$ has no
  roots in~$F$.
\end{lemma}

There are versions of these three lemmas in which the polynomial $T^{p^e}-T$ is replaced by
any separable monic additive polynomial in $F[T]$ of
degree $p^e$ which has all its roots in~$F$ (cf.\ e.g.\ \cite{E}).
However, the statements are simplest in the case of $T^{p^e}-T$,
so we restrict to this case in what follows.

Our next results concern the ramification in $E/F$.

\begin{lemma}  
\label{four}
Suppose $F=k(w)$, where $k$ is perfect (and $w$ is transcendental over $k$).
Let $E=F(v)$, where the minimal polynomial for $v$ over $F$ is
$T^{p^e}-T-z$ (with $z\in F$).  If $E/F$ is unramified over each
finite place of~$F$, then there is some $v'\in E$ for which $v'-v\in F$ and
$E=F(v')$, and where moreover the minimal polynomial for $v'$ over $F$ is
$T^{p^e}-T-z'$ with $z'\in k[w]$.
\end{lemma}

\begin{proof}
For $y\in k(w)$, we may replace $v$ by $v+y$ and $z$ by $z+y^{p^e}+y$
without affecting the hypotheses.  Writing $z=a/b$ with $a,b\in k[w]$
coprime, we may thus assume that $y$ (and hence $z$) has been chosen
to minimize $\deg(b)$.  We will show that $z\in k[w]$.

We begin by showing that $b$ is a $p^e$-th power in $k[w]$.  For this,
    consider any place $Q$ of $k(w)$ which contains $b$ 
    (equivalently, any irreducible $c\in k[w]$ which divides $b$) and let 
    $P$ be a place of $E$ lying over $Q$.
    Let $\rho$ be the (additive) valuation corresponding to~$P$, normalized
    so that $\rho(E)=\bZ$.  Since $P$ is unramified over $Q$, the value
    $\rho(h)$ (for any $h\in k[w]$) is just the multiplicity of $c$ as a
    divisor of $h$.  Thus
    \[\rho(v^{p^e}-v)=\rho(z)=-\rho(b)<0,\]
    so $\rho(v)<0$, whence $p^e\rho(v)=\rho(v^{p^e}-v)=-\rho(b)$.  Since this holds
    for every irreducible factor $c$ of $b$, we conclude that
    $b$ is indeed a $p^e$-th power in $k[w]$; say $b=b_0^{p^e}$.

    For any irreducible factor $c\in k[w]$ of $b$, the residue field
    $k[w]/(c)$ is a finite extension of $k$, hence is perfect; thus,
    $a$ is a $p^e$-th power in this field.  In other words, there is
    some  $a_0\in k[w]$ such that $a_0^{p^e}-a$ is divisible by $c$.
    For $\tilde{v}:=v-a_0/b_0$, we have
    \[ \tilde{v}^{p^e}-\tilde{v}=\frac{a-a_0^{p^e}+a_0b_0^{p^e-1}}{b}.\]
    Since $c$ divides both the numerator and denominator of the
    right-hand expression, we can write $\tilde{v}^{p^e}-\tilde{v}
    = \tilde{a}/\tilde{b}$, where $\tilde{a},\tilde{b}\in k[w]$ 
    and $\deg(\tilde{b})<\deg(b)$.  But this contradicts our
    hypothesis that $\deg(b)$ is minimal.  This implies that $b$ has no
    irreducible factors in $k[w]$, so $b\in k^*$; thus $z\in k[w]$, as desired.
\end{proof}

This result shows that, in the case of interest to us, we may assume
$z\in k[w]$.  Our next two results describe the ramification in extensions
of this sort. 

\begin{lemma}  
\label{five}
Suppose $F=k(w)$, where $k$ is perfect.
Let $E=F(v)$, where the minimal polynomial of
$v$ over $F$ is $h(T):=T^{p^e}-T-z$ (where $z\in k[w]$).  Then $E/k(w)$ is
unramified over each finite place of $k(w)$.  If $k$ is
algebraically closed in $E$ and\/ $\F_{p^e}\subseteq k$, then $E/k(w)$ is 
totally ramified over the infinite place of $k(w)$.
\end{lemma}

\begin{proof} First consider a finite place of $k(w)$, and let $S$ be
  the corresponding valuation ring of $k(w)$; since $h(T)\in S[T]$ 
  and $h'(T)=-1\in S^*$, the place is unramified in the extension $E/k(w)$. 
  Now assume that $k$ is algebraically closed in $E$ and
  $\F_{p^e}\subseteq k$.  By
  Lemma~\ref{two}, $E/k(w)$ is abelian, so each subgroup of $\Gal(E/k(w))$
  is normal.  In particular, this applies to the inertia group
  $I$ of a place $P$ of $E$ lying over the infinite place of 
  $k(w)$, so $I$ does not depend on the choice of $P$.  Thus the 
  fixed field $E^I$ is an unramified Galois extension of 
  $k(w)$ in which $k$ is algebraically closed, 
  so it is a trivial extension; hence $I=\Gal(E/k(w))$,
  so the infinite place of $k(w)$ is totally ramified in $E/k(w)$.
\end{proof}

Now we examine more carefully the ramification over the infinite place, 
by studying the higher ramification groups.  In particular, we consider
the case where there is only one jump in the filtration of ramification
groups over
the infinite place of~$k(w)$; this is what occurs in the situation of
greatest interest to us, namely when $z$ is a power of $w$.  

\begin{prop}  
\label{six}
Suppose $F=k(w)$, where $k$ is perfect and $k\supseteq\F_{p^e}$.
Let $E=F(v)$,
where the minimal polynomial of $v$ over $F$ is 
$T^{p^e}-T-z$ (with $z\in k[w]$).
Assume that no term of $z$ has degree a positive multiple
of~$p^e$, and that $k$ is algebraically closed in $E$.
Let $I_0, I_1,\dots$ be the ramification groups at a place of $E$
lying over the infinite place of $k(w)$.  Then\/ {\em(a)} and\/ {\em(b)} below
are equivalent, and each implies\/ {\em(c)}: 
\renewcommand{\theenumi}{{\alph{enumi}}}
\begin{enumerate}
\item $I_0=I_1=\dots=I_n$ but $I_{n+1}=1$.
\item For every $\zeta\in\F_{p^e}^*$, the polynomial gotten
 {}       from $\zeta z$ by replacing each term $\alpha w^{p^ij}$ by 
          $\alpha^{p^{-i}}w^j$ (for integers $i,j\geq 0$ with $j$ coprime 
          to $p$) has degree~$n$.
\item $n$ is the largest integer coprime to~$p$ which divides 
          the degree of a nonconstant term of $z$.
\end{enumerate}
\end{prop}
\renewcommand{\theenumi}{{\roman{enumi}}}

We may assume that no term of $z$
has degree a positive multiple of $p^e$: we can replace $v$ 
by $\tilde v := v + y$,
where $y\in k[w]$ is chosen so that $\tilde z := 
z+y^{p^e}-y$ has
the desired property (such an element $y$ exists because $k$ is perfect),
and then $F(\tilde v)=F(v)$ and the minimal polynomial of $\tilde v$
over $F$ is $T^{p^e} - T-\tilde z$.  

\begin{proof} By the previous result and Lemma~\ref{two},  
 $I_0=\Gal(E/F)\cong \F_{p^e}$.  Lemma~\ref{two} describes the intermediate 
 fields $F\subset F_0\subseteq E$ for which $[F_0\col F]=p$: they 
 are the fields $F_0=F(v_\zeta)$ where, for $\zeta\in\F_{p^e}^*$, 
 we put $v_\zeta:=(\zeta v)^{p^{e-1}}+(\zeta v)^{p^{e-2}}+\dots+
 (\zeta v)^p+(\zeta v)$.  The minimal polynomial
 for $v_\zeta$ over $F$ is $T^p-T-\zeta z$.  

 First we show the equivalence of (a) and (b).  We begin by relating the
 groups $I_i$ to the corresponding groups for the extensions
 $F(v_\zeta)/F$, via \cite[Prop.~IV.3]{Se} (note that in the statement of
 that result, $e_{L/K}$
 should be replaced by $e_{L/K'}$).
 The stated result 
 implies that (a) is equivalent to the following: for each~$\zeta$, 
 the ramification groups for $F(v_\zeta)/F$ over the infinite place 
 of~$F$ equal $\Gal(F(v_\zeta)/F)$ until the $(n+1)$-th group, which is 
 trivial.
%
%
 Now we just need to compute the ramification groups over infinity
 for the Artin-Schreier extensions $F(v_\zeta)/F$.  We first normalize 
 $v_\zeta$ by replacing it by $\tilde v_\zeta:=v_\zeta+y_\zeta$,
 where $y_\zeta\in k[w]$ is chosen so that no nonzero monomial in
 ${\tilde v_\zeta}^p - \tilde v_\zeta = \zeta z + y_\zeta^p-y_\zeta$ has 
 degree a positive multiple of~$p$; this is possible because $k$
 is perfect.  Then $F(v_\zeta)=F(\tilde v_\zeta)$ and the minimal 
 polynomial for $\tilde v_\zeta$ over $F$ is $T^p-T-\tilde z$, where 
 $\tilde z$ is gotten from $\zeta z$ by replacing each term 
 $\alpha w^{p^i j}$ (where $j$ is coprime to $p$)
 by the term $\alpha^{p^{-i}}w^j$.  The computation of the ramification
 in $F(\tilde v_\zeta)/F$ is classical in this case, since
 $\tilde n:=\deg(\tilde z)$ is coprime to~$p$: the first several ramification
 groups over infinity equal the full Galois group, until the
 $(\tilde n+1)$-th group which is trivial \cite[Prop.~III.7.8]{St}.  This
 proves the equivalence of (a) and (b).  Note that $n=\tilde n$ is
 necessarily coprime to~$p$.

 We now show that (b) implies (c).  So, assume (b); then
 $z$ has a term of degree $np^r$,
 and we have seen that $n$ is coprime to~$p$.  We
 just need to show that $z$ has no term of degree $sp^j$ with 
 $s>n$ and $s$ coprime to~$p$.  For any such $s$, write 
 $z=\sum \alpha_i w^i$ with $\alpha_i\in k$, and
 let $\hat z = \alpha_s w^s + \alpha_{sp} w^{sp} + \alpha_{sp^2} w^{sp^2} +\dots$
 be the sum of the terms of $z$ having degree $sp^j$.  
 Then (b) implies that, for every $\zeta\in\F_{p^e}^*$, we have
 $\zeta \alpha_s + (\zeta \alpha_{sp})^{1/p} + (\zeta \alpha_{sp^2})^{1/p^2} +\dots=0.$
 {}From our assumption that no term
 of $z$ has degree a positive multiple of $p^e$, we see that
 $\deg(\hat z)\leq sp^{e-1}$.  Thus, when we raise the
 previous equation to the power $p^{e-1}$, we get $h(\zeta)=0$,
 where $h(T):=\alpha_s^{p^{e-1}} T^{p^{e-1}} + \alpha_{sp}^{p^{e-2}} T^{p^{e-2}}
 + \dots$ is a polynomial in $k[T]$.  Since $h(T)$ vanishes on $\F_{p^e}$,
 but $\deg(h)\leq p^{e-1}$, it follows that $h=0$, so each 
 $\alpha_{sp^j}=0$.  This concludes the proof. 
\end{proof}

We now show when certain data determines the field
$k(v,w)$, where $v^{p^e}-v=w^n$.  The data come from the extension
$k(v,w)/k(w^r)$ for some $r$; we 
are given the ramification groups over infinity for this extension, 
and also we are given that $k(v,w)/k(w)$ is unramified over 
finite places.  We will see that this data uniquely determines
$k(v,w)$ if and only if $p,n,r$ satisfy a certain arithmetic
condition.  

\renewcommand{\theenumi}{{\alph{enumi}}}
\begin{prop}  
\label{msri}
Let $E\supseteq k(\tilde w)\supseteq F$ satisfy
$[E\col k(\tilde w)]=p^e$ and $[k(\tilde w)\col F]=r$, where $k$ is a perfect 
field containing\/ $\F_{p^e}$ and $\tilde w$ is transcendental over~$k$.
Assume that $E/F$ is Galois and is totally ramified over some
degree-one place of~$F$, and 
that the sequence of ramification groups over this
place has jumps only after the 0-th and $n$-th groups (where $p\nmid n$).
Also assume that $E/k(\tilde w)$ is ramified over only one place, and that
$r/\!\gcd(n,r)$ divides $p^e-1$.  Finally, assume that either $k=\F_{p^e}$ or
$p^e$ is the least power of $p$ congruent to 1 modulo~$r/\!\gcd(n,r)$.
Then $E=k(v,w)$ where
\begin{enumerate}
\item $k(w)=k(\tilde w)$ and $F=k(w^r)$;
\item $z:=v^{p^e}-v$ lies in $k[w]$ and has no nonconstant terms
  of degree divisible by~$p^e$;
\item $n$ is the largest integer coprime to $p$ which
   divides the degree of a nonconstant term of $z$;
\item $E/k(w)$ is Galois, and the map
  $\sigma\mapsto\sigma(v)-v$
  induces an isomorphism\/ $\Gal(E/k(w))\to \F_{p^e}$; for $\alpha\in\F_{p^e}$,
  let $\sigma_\alpha$ be the corresponding element of\/ $\Gal(E/k(w))$;
\item for $\tau\in\Gal(E/F)$, put $\zeta:=\tau(v)/v$; then
   $\zeta\in\F_{p^e}^*$ and $\tau\sigma_\alpha\tau^{-1}=
  \sigma_{\alpha\zeta^{-n}}$;
\item every term of $z$ has degree congruent to $n$ mod~$r$.
\end{enumerate}
\end{prop}
\renewcommand{\theenumi}{{\roman{enumi}}}

The proof relies on various results about ramification groups.  The 
standard reference for these is \cite[Ch.~IV]{Se}; we recall
the facts we will need.  Given 
a Dedekind domain $S$ with field of fractions $F$, let $R$ be its 
integral closure in a finite Galois extension $E$ of $F$, with Galois 
group $G$.  Let $P$ be
a prime ideal of $R$, and put $Q=P\cap S$.  The decomposition group $D$
of $P$ is the subgroup of $G$ consisting of elements $\sigma\in G$ with
$\sigma(P)=P$.  When the extension of residue fields $(R/P)/(S/Q)$ is 
separable, the $i$-th ramification group $I_i$ of $G$ relative to $P$ 
(for $i\geq 0$) is
defined to be the set of $\sigma\in G$ which act trivially on $R/P^{i+1}$.
The $I_i$ form a decreasing sequence of normal 
subgroups of $D$, and $I_i=1$ for $i$ sufficiently large.
Here $I_0$ is the inertia group of $P$.
By a `jump' in the sequence of ramification groups, we mean an
integer $i\ge 0$ for which $I_i\neq I_{i+1}$.
Especially important for our purposes are results of \cite[IV.2]{Se},
which we now state (our statements differ slightly from those of \cite{Se}).
We denote the residue field $R/P$ by~$\ell$.

\begin{lemma}  
\label{serre}
The map $\theta_0\colon I_0/I_1\to\Aut_\ell(P/P^2)$
 given by $\theta_0(\tau)\colon\pi\mapsto\tau(\pi)$ is an injective
 homomorphism.  For $i\geq 1$, the map 
 $\theta_i(\sigma)\colon\pi\mapsto\sigma(\pi)-\pi$
 induces an injective homomorphism
 $\theta_i\colon I_i/I_{i+1} \to \Hom_\ell(P/P^2,P^{i+1}/P^{i+2})$.
 For $\tau\in I_0$ and $\sigma\in
 I_i/I_{i+1}$, we have $\theta_i(\tau\sigma\tau^{-1})=\theta_0(\tau)^i
 \theta_i(\sigma)$.
\end{lemma}

Here $P/P^2$ is a one-dimensional $\ell$-vector space, so
$\Aut_\ell(P/P^2)\cong\ell^*$.  Likewise $P^{i+1}/P^{i+2}$ is a 
one-dimensional $\ell$-vector space, so $P^i/P^{i+1}$ is isomorphic to
$\Hom_\ell(P/P^2,P^{i+1}/P^{i+2})$ via the map taking $\psi$ to
$\pi\mapsto\psi\pi$.  Thus the right side of the final equation in
the lemma makes sense, since it is just the action of $\ell^*$ on
the $\ell$-vector space $P^i/P^{i+1}$.  Finally, note that the
final equation simply amounts to the natural action 
of $P/P^2$ on its $i$-th tensor power, which explains the $i$-th power
in that equation.

\begin{proof}[Proof of Proposition~\ref{msri}] 
Let $Q$ be the degree-one
place of $F$ over which $E/F$ is totally ramified, and let $P$ be
the place of $E$ lying over $Q$.  Since $P/Q$ is totally ramified,
the inertia group $I_0$ equals $G:=\Gal(E/F)$.  Since in addition
$Q$ has degree one, the constant fields of $F$, $k(\tilde w)$, and
$E$ must all be the same, so they are all~$k$.  
By replacing $\tilde w$ by $1/(\tilde w -\alpha)$ if necessary (with $\alpha\in
k$), we may
assume that $P$ lies over the infinite place of $k(\tilde w)$.
Next, $I_1$ is the unique Sylow $p$-subgroup of $I_0$, 
so the fixed field $E^{I_1}$ must equal $k(\tilde w)$ (since
$[E\col k(\tilde w)]=p^e$).  Hence $\Gal(k(\tilde w)/F)\cong I_0/I_1$ 
is cyclic of order~$r$;
here $\langle\nu\rangle:=\Gal(k(\tilde w)/F)\subset\Aut_k(k(\tilde w))$, and the 
latter group is isomorphic to $\PGL(k)$ where $\textmatrix{\alpha&\beta\\ \gamma&\delta}$ 
corresponds to the $k$-automorphism of $k(\tilde w)$ sending 
$\tilde w\mapsto (\alpha\tilde w+\gamma)/(\beta\tilde w+\delta)$.  
Since $k(\tilde w)/F$ is totally ramified under the infinite place,
$\nu$ generates the decomposition group under this place, so
$\nu(\tilde w)=\alpha{\tilde w}+\gamma$; but $\nu$ has order $r$ which is coprime to
$p$, so either $\nu$ is the identity (and we put $w=\tilde w$) or
$\alpha\neq 1$, in which case we put $w=\tilde w+\gamma/(\alpha-1)$.  In any
case, there exists $w\in k[\tilde w]$ such that
$k(w)=k(\tilde w)$ and $\nu(w)=\alpha w$ with $\alpha\in k$ a primitive $r$-th
root of unity, so $F=k(w)^{\nu}=k(w^r)$.
This proves (a); note that $E/k(w)$ is totally ramified over infinity.

Let $\tau\in\Gal(E/F)$ map to a generator of $\Gal(k(w)/F)$;
then $\tau(w)=\zeta w$ where $\zeta$ is a primitive $r$-th root of
unity.  Since $P$ is totally ramified over the 
infinite place of $k(w)$, and $1/w$ is
a uniformizer for the latter place, we have
$1/w\in P^{p^e}\setminus P^{p^e+1}$.  If $A$ denotes the valuation
ring of $E$ corresponding to $P$, and $\pi$ is any uniformizer of $P$,
then $1/w=y\pi^{p^e}$ for some $y\in A^*$; since $\tau(y)-y\in P$, 
\[\tau(1/w)=
\tau(y)\tau(\pi)^{p^e}\equiv y\tau(\pi)^{p^e}\pmod{P^{p^e+1}}.\]
But $\tau(1/w)=1/(\zeta w) = \zeta^{-1}y\pi^{p^e}$, so 
$\tau(\pi)^{p^e}\equiv \zeta^{-1}\pi^{p^e}\pmod{P^{p^e+1}}$,
whence $\theta_0(\tau)=\zeta^{-1/p^e}\in k^*$.
Next, $\theta_n$ induces an injective homomorphism
$I_1\to P^n/P^{n+1}$, and $\theta_n(\tau\sigma\tau^{-1})=\zeta^{-n}
\theta_n(\sigma)$ for $\sigma\in I_1$.  Iterating, 
$\theta_n(\tau^i\sigma\tau^{-i})=\zeta^{-in}\theta_n(\sigma)$,
so the image of $\theta_n$ contains all sums of elements
$\zeta^{-in}\theta_n(\sigma)$; in other words, the image of $\theta_n$ 
contains $\F_p(\zeta^{-n})\theta_n(\sigma)$, which equals 
$\F_{p^e}\theta_n(\sigma)$ by hypothesis.  Choose $\sigma$ to be any
nonidentity element of $I_1$; then the image of $\theta_n$ is
precisely $\F_{p^e}\theta_n(\sigma)$, and we have an isomorphism
{}from $\F_{p^e}$ to the image of $\theta_n$ via $\alpha\mapsto
\alpha\theta_n(\sigma)$ (unlike our previous isomorphisms, this
one is not canonical).

Now we use basic Galois cohomology to pick the element~$v$.
A reference is \cite[Chs.~VII and X]{Se}.
We have a homomorphism $\rho\colon I_1\to\F_{p^e}\subset E$, which is
a 1-cocycle for the $I_1$-module $E$.  Since
$H^1(I_1,E)=0$, the cocycle $\rho$ is a coboundary, so there exists
$\tilde v\in E$ such that, for each $\sigma\in I_1$, we have
$\rho(\sigma)=\sigma(\tilde v)-\tilde v$.  It follows that the map
$\sigma\mapsto\sigma(\tilde v)-\tilde v$ is an isomorphism 
$I_1\to\F_{p^e}$;
we denote by $\sigma_\alpha$ the preimage of $\alpha\in\F_{p^e}$ under
this map.  In particular, $\tilde{v}$ has $p^e$ conjugates under
$I_1$, so indeed $E=k(\tilde v, w)$.
Also we now know the shape of the minimal polynomial for $\tilde v$ over 
$k(w)$: it is $\prod_{\alpha\in\F_{p^e}}(T-(\tilde v+\alpha))=
(T-\tilde v)^{p^e}-(T-\tilde v)= T^{p^e}-T-\tilde z$, where 
$\tilde z:=\tilde v^{p^e}-\tilde v$ lies in $k(w)$.

By Lemma~\ref{four} and the remark following
Proposition~\ref{six},  there exists 
$v'\in E$ such that $v'-\tilde v\in k(w)$ and $z':={v'}^{p^e}-v'$ 
lies in $k[w]$ but $z'$ has no terms of degree a positive multiple 
of $p^e$.  Note that $v'$ is only determined up to addition by an 
element of~$k$; eventually we will
specify the choice of this element to determine $v$.
Regardless of the choice, we have $k(v',w)=k(\tilde v,w)$ and 
$\sigma_\alpha(v')-v'=\alpha$; (b) and (d) follow at once, and
(c) then follows from Proposition~\ref{six}.  

Let $\tau\in\Gal(E/F)$ map $w\mapsto\zeta w$, where $\zeta$ has
order~$r$.  Since $\theta_n(\tau^i\sigma_\alpha\tau^{-i})=
\zeta^{-in}\theta_n(\sigma_\alpha)$, injectivity of $\theta_n$ yields 
$\tau^i\sigma_\alpha\tau^{-i}=\sigma_{\alpha\zeta^{-in}}$;
since every element of $\Gal(E/F)$ has the form 
$\tau^i\sigma_{\tilde\alpha}$, this proves~(e).  

Finally, we consider (f), and specify the choice of $v\in v'+k$.  
For $\tau$ as in the previous paragraph, we have
\[\alpha\zeta^{-n}=\sigma_{\alpha\zeta^{-n}}(v')-v'
    = \tau\sigma_\alpha\tau^{-1}(v')-v'
    = \tau(\sigma_\alpha\tau^{-1}(v')-\tau^{-1}(v'));\]
since $\tau$ fixes $\alpha\zeta^{-n}$, it follows that
$\alpha\zeta^{-n}=\sigma_\alpha\tau^{-1}(v')-\tau^{-1}(v')$.  But also 
$\alpha\zeta^{-n}=\sigma_\alpha (\zeta^{-n}v')-\zeta^{-n}v'$,
so, subtracting, we see that $y:=\tau^{-1}(v')-\zeta^{-n}v'$ is
fixed by each $\sigma_\alpha$;  hence $y\in k(w)$.
Since $y^{p^e}-y=\tau^{-1}(z')-\zeta^{-n}z'$ lies in $k[w]$
and has no term
of degree a positive multiple of $p^e$, in fact $y$ must lie
in~$k$.  For $\beta\in k$, we have 
\[\tau^{-1}(v'+\beta)=\zeta^{-n}(v'+\beta)+\beta(1-\zeta^{-n})+y;\]
if $\zeta^n\neq 1$, there is a unique choice of $\beta$ for
which $v:=v'+\beta$ satisfies $\tau^{-1}(v)=\zeta^{-n}v$.
If $\zeta^n=1$, then by replacing $\tau$ by its $p$-th power 
(and choosing $v'=v$) we may assume that $\tau^{-1}(v)=v$.
At last, for $z:=v^{p^e}-v\in k[w]$, we have 
\[z\tau^{-1}(z)=\tau^{-1}(v^{p^e}-v)=\zeta^{-n}(v^{p^e}-v)
     = \zeta^{-n}z,\]
so all terms of~$z$ have degree congruent to $n$ mod~$r$.
This completes the proof. 
\end{proof}

\begin{cor}  
\label{coro}
Under the hypotheses of Proposition~\ref{msri},
suppose the integers $p,n,r$ satisfy
\[\hbox{if}\quad n',i\geq0\quad\hbox{and}\quad n'\equiv p^in\!\pmod{r},%
\quad\hbox{then}\quad  n'\geq n. \leqno{(*)}\]
Then conclusions (a)--(f) of Proposition~\ref{msri} hold if
we require in addition that $z=\gamma w^n$ with $\gamma\in k^*$.
\end{cor}

\begin{proof}
Let $v,w,z$ be as in the conclusion of the proposition, so $z\in
k[w]$.  Then the degree of any nonconstant term of $z$ is
congruent to $n$ mod~$r$ (by (f)) and has the shape $n'p^j$
with $0\leq n'\leq n$ (by (c)) and $0\leq j<e$ (by (b)), so 
($*$) implies $n'=n$.  If $p^e$ is the least power of $p$ congruent to 1 
modulo~$r/\gcd(n,r)$, then ($*$) implies $j=0$ and we are done.
If this condition does not hold then, by the hypothesis of
the proposition, $k=\F_{p^e}$.  In this case, write
$z=\sum_{j=0}^{e-1} \alpha_j w^{np^j}$, where each $\alpha_j\in\F_{p^e}$
and some $\alpha_j\neq 0$ (by (c)).  Let $\hat v$ be an element of an
extension of $E$ satisfying ${\hat v}^{p^e}-\hat v=w^n$; then
$\tilde{y}:=\sum_{j=0}^{e-1}\alpha_j {\hat v}^{p^j}$ satisfies
$\tilde{y}^{p^e}-\tilde{y}=z$, so $\tilde{y}-y\in\F_{p^e}$.
Since $[k(\hat v,w)\col k(w)]=p^e$ (by Lemma~\ref{three})
and $k(\hat v,w)\supseteq k(v,w)$, 
we have $k(v,w)=k(\hat v,w)$, and the result follows.
\end{proof}

We now present values of the parameters for which ($*$) is satisfied.

\begin{lemma}  
\label{numerology}
The criterion {\em($*$)} is satisfied if $p^e\geq 4$ and
either of the following hold:
\begin{enumerate}
\item $p^e\equiv 3\!\pmod 4$, 
    $\,\,n\mid ((p^e+1)/4)$, $\,\,r=(p^e-1)/2$.
\item $p^e\not\equiv 3\!\pmod 4$, $\,\,n\mid p^e+1$, 
  $\,\,n<p^e+1$, $\,\,r=p^e-1$.
\end{enumerate}
\end{lemma}

\begin{proof} It suffices to prove (i) in case $n=(p^e+1)/4$, since we 
can reduce the
general case to this one by multiplying ($*$) by $(p^e+1)/(4n)$.
Now, \[n\cdot p=n+\frac{p^e+1}{2}\cdot\frac{p-1}{2}
       \equiv n+\frac{p-1}{2}\pmod{r},\]
so \[n\cdot p^i\equiv n+\frac{p-1}{2}\cdot
     (1+p+\dots+p^{i-1})=n+\frac{p^i-1}{2}\pmod{r}.\]
For $0\leq i<e$ we have \[\frac{p^e+1}{4}+\frac{p^i-1}{2}\leq\frac{p^e+1}{4}+
  \frac{p^e/3-1}{2}=\frac{5p^e-3}{12}<r,\]
so $n\,+\,(p^i-1)/2$ is the least nonnegative residue of
$np^i$ modulo~$r$.  In particular, the least such number is~$n$;
since $p^e\equiv 1\!\pmod{r}$, this proves~(i).
  
Now consider (ii).  We assume $n>1$, since the result is clear for $n=1$.
Assume $0\leq i<e$.
Multiplying ($*$) by $(p^e+1)/n$ gives
\[\frac{p^e+1}{n}\cdot n'\equiv p^i(p^e+1)\pmod{r}.\]
The right side is congruent to $p^{i+1}$, and it follows that its
least nonnegative reside modulo~$r$ is itself (unless $p=2$ and
$i=e-1$ in which case it is 1); but if
$n'<n$ then the left side is already reduced modulo~$r$, so 
the only possibility is $p$ odd, $(p^e+1)/n=2$, and $n'=p^i$.  In
that case $n\equiv 1\!\pmod{r}$, contradiction.
\end{proof} 

The specific values in this lemma are the ones we will use elsewhere
in this paper.  For completeness, we now present necessary and
sufficient conditions for the fields $k(v,w)$ (with $v^q-v=\gamma w^n$)
to be determined by their ramification over $k(w^r)$.  We first prove
Lemma~\ref{lem}, which describes this ramification; this lemma has
been known for many years.

\begin{proof}[Proof of Lemma~\ref{lem}]
By Lemma~\ref{three}, $T^q-T-\gamma w^n$ is irreducible in $\bar{k}(w)[T]$,
so $[k(v,w)\col k(w)]=q$ and $\bar{k}(w)\cap k(v,w)=k(w)$.  
It follows that $k(v,w)/k(w^r)$ has degree $qr$; suppose this
extension is Galois with group $G$.  Since it is separable, $r$ is
coprime to $q$.  Since it is normal, there must be elements of $G$
which map $w$ to any of its conjugates over $k(w^r)$; in particular,
there exists $\sigma\in G$ with $\sigma(w)=\zeta w$, where $\zeta$
is a primitive $r$-th root of unity (and $\zeta\in\bar{k}\cap
k(v,w)=k$).  Write $\sigma(v)=\sum_{i=0}^{q-1} y_i v^i$ with
$y_i\in k(w)$; then 
\[\gamma\zeta^n v^n = \sigma(v)^q-\sigma(v)=\sum_{i=0}^{q-1}
  (y_i^q(v+\gamma w^n)^i-y_iv^i).\]
By considering the terms of degree $i$ in $v$, for $i=q-1,
q-2,\dots,0$ successively, we find that $y_i=0$ for $i>1$, and also
$y_0\in\F_q$ and $\zeta^n=y_1\in\F_q$.  This last statement may
be restated as: $r/\gcd(n,r)$ divides $q-1$.

Conversely, assume that $k$ contains a primitive $r$-th root of unity
and $r/\gcd(n,r)$ divides $q-1$.  Let $\zeta$ be an $r$-th root of unity
and let $\alpha\in\F_q$; then $\zeta^n\in\F_q$.  Then there is a
$k$-automorphism of $k(v,w)$ mapping $w\mapsto\zeta w$ and
$v\mapsto \zeta^n v+\alpha$ (since these equations define an
automorphism of $k(w)[X]$, and they preserve the ideal generated
by $X^q-X-\gamma w^n$, so they define an automorphism of the quotient ring
which is $k(v,w)$).  These automorphisms form a group of order $qr$,
and they all fix $w^r$ (where $[k(v,w)\col k(w^r)]=qr$), 
so $k(v,w)/k(w^r)$ is Galois.  The ramification in this extension
is as in ($\dagger$), by Lemma~\ref{five}
and Proposition~\ref{six}.
\end{proof}

We conclude this section by proving Theorem~\ref{unq}.
\begin{proof}[Proof of Theorem~\ref{unq}]
The first part of the theorem follows from Corollary~\ref{coro}
once we observe that the subfield $E^{I_1}$ has the
form $k(\tilde w)$; for this, note that $E^{I_1}/k(t)$
is cyclic of degree $r$ with only two branch points
($0$ and infinity), both of which have degree one and are totally (and tamely)
ramified, so $E^{I_1}$ is a genus zero function field over $k$ having
a degree one place, whence it indeed has the form $k(\tilde w)$.

Assume that (i) is violated, and let $0\leq n'<n$ satisfy 
$np^i\equiv n'\!\pmod{r}$ for some $i\geq 0$; since
$np^e\equiv n\!\pmod{r}$, we may assume $0\leq i<e$.
Let $v$ and $w$ be transcendentals over $k$ satisfying
$v^{p^e}-v=w^{np^i}+w^{n'}$.
By Lemma~\ref{three}, $[\bar{k}(v,w)\col \bar{k}(w)]=p^e$.
For any $r$-th root of unity $\zeta$ and any $\alpha\in\F_{p^e}$,
there is a $k$-automorphism of $k(v,w)$ mapping 
$w\mapsto\zeta w$ and $v\mapsto \zeta^{n'}v+\alpha$.  These 
automorphisms form a group of order $p^er$, and their fixed field is
$k(w^r)$, so $k(v,w)/k(w^r)$ is Galois.  The ramification in
this extension is as in ($\dagger$), by Lemma~\ref{five} and
Proposition~\ref{six}.  Now suppose that $k(v,w)=k(\hat{v},\hat{w})$
where $k(w^r)=k(\hat{w}^r)$ and $\hat{v}^{p^e}-\hat{v}=\gamma\hat{w}^n$ for some
$\gamma\in k^*$.
Then we would have $w=\beta \hat{w}$ with $\beta\in k^*$ and
$v=\sum_{i=0}^{p^e-1}y_i \hat{v}^i$ with $y_i\in k(\hat{w})$, so
\[
(\beta \hat{w})^{np^i}+(\beta \hat{w})^{n'}=v^{p^e}-v 
=\sum_{i=0}^{p^e-1}\left(y_i^{p^e}(\hat{v}+\gamma\hat{w}^n)^i-y_i\hat{v}^i\right),
\]
and we get a contradiction by considering successively the
terms of degree $i$ in $\hat{v}$ for $i=p^e-1,p^e-2,\dots,0$.

Now assume (ii) is violated.  Pick an integer $1\leq i<e$ such that
$p^i\equiv 1\!\pmod{r/\gcd(n,r)}$, and let $v,w$ be transcendentals
over $k$ satisfying $v^{p^e}-v=w^n+\beta w^{np^i}$, where 
$\beta\in k\setminus\F_q$.  The above argument shows that
$k(v,w)/k(w^r)$ is Galois with ramification as in ($\dagger$),
and that this extension cannot be written in the form
$k(\hat{v},\hat{w})/k(\hat{w}^r)$ for any $\hat{v},\hat{w}\in k(v,w)$ with
$\hat{v}^{p^e}-\hat{v}=\gamma\hat{w}^n$ and $\gamma\in k^*$.
\end{proof}

\section{Producing polynomials}

In this section we use the results from the previous two
sections to prove refined versions of the results stated in
the introduction.  We first describe the polynomials $f(X)$ of degree
$p^e(p^e-1)/2$ (over a field $k$ of characteristic $p>0$)
whose arithmetic monodromy group has a transitive
normal subgroup isomorphic to $\PSL(p^e)$, assuming that either
the Galois closure of the extension $k(x)/k(f(x))$ 
does not have genus $p^e(p^e-1)/2$, or this extension has no finite
branch points.  We then apply this result to produce all
indecomposable polynomials $f$ (over any field $k$) such that 
$\deg(f)$ is not a power of $\charp k$ and either $f$ is
exceptional over $k$ or $f$ decomposes over an extension
of $k$.  We begin with the following refinement of Theorem~\ref{genpol}:

\renewcommand{\theenumi}{{\arabic{enumi}}}
\begin{thm}
\label{genpol2}
Let $k$ be a field of characteristic $p>0$,
let $d=(q^2-q)/2$ for some power $q=p^e$, and let 
$f(X)\in k[X]\setminus k[X^p]$ have degree $d$.
Then the following are equivalent:
\begin{enumerate}
\item $\Gal(f(X)-u,k(u))$ has a transitive normal subgroup isomorphic to\/
$\PSL(q)$, and the extension $k(x)/k(f(x))$ either has no finite
branch points, or has Galois closure of genus $\neq(q^2-q)/2$.
\item There exist linear polynomials $\ell_1,\ell_2\in \bar{k}[X]$
   such that the composition $\ell_1\circ f\circ\ell_2$ is one
   of the following polynomials or one of the exceptions in Table~B below:
\[X(X^m+1)^{(q+1)/(2m)}
 \left(\frac{(X^m+1)^{(q-1)/2}-1}{X^m}\right)^{(q+1)/m}\]
with $q$ odd and $m$ a divisor of $(q+1)/2$; or
\[X^{-q}\biggl(\sum_{i=0}^{e-1}X^{m2^i}\biggr)^{(q+1)/m}\]
with $q$ even and $m$ a divisor of $q+1$.
\end{enumerate}
In these cases, $\Gal(f(X)-u,\bar{k}(u))$ is either $\PSL(q)$ (if $m$ is even)
or $\PGL(q)$ (if $m$ is odd).  Moreover, the Galois closure of
$\bar{k}(x)/\bar{k}(f(x))$ is $\bar{k}(v,w)$ where $v^q-v=w^{m/\!\gcd(m,2)}$.
\end{thm}
\renewcommand{\theenumi}{{\roman{enumi}}}
{\tiny{
\begin{align*}
f_1 = \,\,& X^3 (X^2+3)^2 (X^8-X^6+2X^4-X^2-3)^6 \\  
f_2 = \,\,& X^3 (X^4+3) (X^{16}-X^{12}+2X^8-X^4-3)^3 \\  
f_3 = \,\,& X^3 (X+3)^4 (X^4-X^3+2X^2-X-3)^{12} \\     
f_4 = \,\,& (X^3-3) X^4 (5+6X^3-X^6-2X^9+X^{12})^4 \\  
f_5 = \,\,& (X-11)^3 (X-10)^4 (X+3)^6 
                (X^{20}-5X^{19}-10X^{18}+7X^{17}+X^{16}+5X^{15}-9X^{14}-10X^{13}-
8X^{12}\\&+11X^{11}
                 +8X^9+10X^8+X^7-6X^6-8X^5+6X^4+12X^3+6X^2-X-7)^{12} \\
f_6 = \,\,& X^3 (X^2+1)^2 (X^2-9)^3 (-3-8X^2-8X^4+11X^6+10X^8-10X^{12}-2X^{14}+7X^{16}- 
            6X^{18}+5X^{20}\\&-4X^{24}+8X^{26}+2X^{28}+2X^{30}+6X^{32}+11X^{34}-7X^{36}+
8X^{38}+X^{40})^6\\
f_7 = \,\,& (X^3-1) X^4 (X^3-10)^2
(-2-X^3+11X^6-X^9+X^{12}-9X^{15}-5X^{18}+X^{24}+2X^{27}- 
             6X^{30}+2X^{33}\\&-6X^{36}-3X^{39}+6X^{42}+9X^{45}+4X^{48}-3X^{51}+8X^{54}+
11X^{57}+X^{60})^4
			 \\
f_8 = \,\,& (X-41)^6  (X-6)^{10} (X-9)^{15}
(19X + 46X^2 + 21X^3 + 2X^4 + 2X^5 + 16X^6 + 53X^7 + 8X^8 + 42X^9 + 
             22X^{10} \\&+ 14X^{11} + X^{12} + 4X^{13} + 12X^{14} +
			 33X^{15} + 41X^{16} +  50X^{17} + 
             27X^{18} + 37X^{19} + 42X^{20} + 8X^{21} + 16X^{22}\\& +
			 53X^{23} +  28X^{24} + 
             9X^{25} + 56X^{26} + 39X^{27} + 42X^{28} + 13X^{29} + 14X^{30} + 28X^{31} + 
             25X^{32} + 26X^{33}\\& + 43X^{34} + 34X^{35} + 10X^{36} + 17X^{37} + 58X^{38} + 
             25X^{39} + 48X^{40} + 14X^{41} + 15X^{42} + 53X^{43} + 39X^{44} \\&+ 58X^{45} + 
             48X^{46} + 5X^{47} + 8X^{48} + 9X^{49} + 9X^{50} + 9X^{51} + 27X^{52} + 4X^{53} + 
             13X^{54} + 56X^{55} + X^{56})^{30}   \\
f_9 = \,\,& (X^2-32)^3 X^{15} (X^2+3)^5 
(43 + 13X^2 + 15X^4 + 2X^6 + 57X^8 + 57X^{10} + 15X^{12} + 43X^{14} + 
     35X^{16} + 53X^{18} \\&+ 35X^{20} + 37X^{22} + 51X^{24} + 5X^{26} + 44X^{28} + 
     6X^{30} + 3X^{32} + 28X^{34} + 44X^{36} + 26X^{38} + 57X^{40} \\&+40X^{42} + 
     41X^{44} + 11X^{46} + 28X^{48} + 57X^{50} + 3X^{52} + 56X^{54} +
	 6X^{56} + 
     30X^{58} + 16X^{60} + 26X^{62} \\&+ 15X^{64} + 36X^{66} + 19X^{68}
	 + X^{70} +  7X^{72} + 
     53X^{74} + X^{76} + X^{78} + 34X^{80} + 32X^{82} + 16X^{84} + 28X^{86} \\&+ 12X^{88} + 
     15X^{90} + 53X^{92} + 20X^{94} + 8X^{96} + 8X^{98} + 14X^{100} + 53X^{102} + 
     54X^{104} + 24X^{106} + 17X^{108} \\&+ 29X^{110} + X^{112})^{15} \\ 
f_{10} = \,\,& (X^3-35)^2 X^{10} (X^3-3)^5
(27 + 45X^3 + 30X^6 + 26X^9 + 41X^{12} + 24X^{15} + 16X^{18} + 43X^{21} + 
     7X^{24} \\&+ 39X^{27} + 24X^{33} + X^{36} + 32X^{39} + 47X^{42} + 37X^{45} + 38X^{48} + 
     18X^{51} + 16X^{54} + 7X^{60} + 24X^{63} \\&+ 48X^{66} + 8X^{69} + 54X^{72} + 
     56X^{75} + 36X^{78} + X^{81} + 33X^{84} + 35X^{87} + 31X^{90} + 34X^{93} + 
     19X^{96} \\&+ 17X^{99} + 29X^{102} + 25X^{105} + 16X^{108} + 17X^{111} + 2X^{114} + 
     8X^{117} + 46X^{120} + 53X^{123} + 54X^{126} \\&+ 15X^{129} + 24X^{132} + 2X^{135} + 
     49X^{138} + 22X^{141} + 36X^{144} + 36X^{147} + 51X^{150} +
	 15X^{153} + 9X^{156}  \\&+ 
     32X^{159} + 6X^{162} + 38X^{165} + X^{168})^{10} \\ 
f_{11} = \,\,& X^6 (X^5+35)^2 (X^5+32)^3
(54 + 12X^5 + 49X^{10} + 33X^{15} + 13X^{20} + 31X^{25} + 13X^{30} + 
      32X^{35} + 6X^{40} \\&+ 10X^{45} + 43X^{50} + 11X^{60} + 54X^{65}
	  + 40X^{70} +  49X^{75} + 
      X^{80} + 13X^{85} + 37X^{90} + 49X^{95} + 40X^{100} \\&+ 10X^{105}
	  + 43X^{110} +  2X^{120} + 
      24X^{125} + 54X^{130} + 46X^{135} + 8X^{140} + 33X^{145} + 35X^{150} + 23X^{155} \\&+ 
      2X^{160} + 57X^{165} + 15X^{170} + 30X^{180} + 32X^{185} + 39X^{190} + 50X^{195} + 
      50X^{200} + 36X^{205} + 55X^{210} \\&+ 15X^{215} + 30X^{220} + 35X^{225} + 3X^{230} + 
      53X^{235} + 37X^{240} + 52X^{245} + 31X^{250} + 6X^{255} + 35X^{260} \\&+ 37X^{265} + 
      30X^{270} + 51X^{275} + X^{280})^6 
\end{align*}}}
\noindent ${}$\qquad\qquad\qquad\qquad\qquad\qquad {{Table~B}}
\newline
\noindent

The first four polynomials in Table~B correspond to $q=11$, the next three
correspond to $q=23$, and the final four correspond to $q=59$.  The
geometric monodromy group $G$ of each polynomial in Table~B is $\PSL(q)$, except
for $f_3$ and $f_4$ where $G\cong\PGL(q)$.  The 
point-stabilizer $G_1$ of $G$ is $A_4$ for the first two polynomials, $S_4$ 
for the next five polynomials, and $A_5$ for the final four polynomials.
The Galois closure of $\bar{k}(x)/\bar{k}(f_i(x))$ has the form $\bar{k}(v,w)$ with
$v^q-v=w^m$ and $m=1$ (for $i=1,3,5,8$) or $m=2$ (for $i=2,6,9$)
or $m=3$ (for $i=4,7,10$) or $m=5$ (for $i=11$).

We postpone the proof of Theorem~\ref{genpol2} until Section~\ref{subsec}.
Here is an outline of the strategy.
Let $f(X)\in k[X]$ satisfy the hypotheses of the theorem, and
assume that (1) holds; then Theorem~\ref{groups} lists the
possibilities for the ramification in $E/\bar{k}(f(x))$,
where $E$ denotes the Galois closure of $\bar{k}(x)/\bar{k}(f(x))$.
For each ramification possibility, let $G$ and $G_1$ be the
corresponding possibilities for the geometric monodromy group and a
one-point stabilizer.
Let $I$ be the inertia group of a place of $E$
lying over the infinite place of $\bar{k}(f(x))$.  First we show that
$E^I$ has genus zero, and that the ramification in 
$E/E^I$ is described by ($\dagger$) (cf.\ Lemma~\ref{lem}).
Then Theorem~\ref{unq} and Lemma~\ref{numerology} imply that  $E=\bar{k}(v,w)$ where
$v^{p^e}-v=w^n$.  The automorphism
group $\Aut_{\overline{k}}(E)$ was determined by Stichtenoth~\cite{St2};
in our cases, it has a unique conjugacy class of subgroups
isomorphic to $G$.  We compute the subfield $\bar{k}(\hat u)$ of $E$ invariant
under one such subgroup, and also the subfield $\bar{k}(\hat x)$ invariant
under the one-point stabilizer $G_1$.  Finally, we compute
the rational function $\hat f$ for which $\hat u=\hat{f}(\hat x)$; for
suitably chosen $\hat u$ and $\hat x$, the rational function $\hat f$ will in fact
be a polynomial $f$ satisfying (1).  Moreover, this construction produces all polynomials
satisfying (1).

\subsection{Consequences of Theorem 4.1}

We now use Theorem~\ref{genpol2} to deduce the remaining
results mentioned in the introduction.
First consider Theorem~\ref{fgsthm}.  Our next results classify the
examples in case (iii) of Theorem~\ref{fgsthm}, assuming that item (1)
of Theorem~\ref{genpol2} holds.  (We will remove this assumption in
\cite{GRZ}, which leads to a new family of examples in characteristic $2$.)
We classify these polynomials up to \emph{equivalence},
where we say that $b,c\in k[X]$ are equivalent if $b=\ell_1\circ c
\circ \ell_2$ for some linear polynomials $\ell_1,\ell_2
\in k[X]$; trivially, this equivalence relation preserves
indecomposability, exceptionality, and both the arithmetic and
geometric monodromy groups.

\begin{thm}
\label{char3thm}
Let $k$ be a field of characteristic $3$, and let $q=3^e$ with $e>1$ odd.
The following are equivalent:
\begin{enumerate}
\item there exists an indecomposable exceptional polynomial $f\in k[X]$
of degree $q(q-1)/2$ for which\/ $\PSL(q)$ is a transitive normal subgroup of
$\Gal(f(X)-u,k(u))$;
\item $k\cap\F_{q}=\F_3$ and $k$ contains non-square elements.
\end{enumerate}
Moreover, for any $f$ as in {\rm (i)}, the Galois closure of
$\bar{k}(x)/\bar{k}(f(x))$ is $\bar{k}(v,w)$ where $v^{q}-v=w^n$
and $n$ divides $(q+1)/4$; here $n$ is uniquely determined by $f$, and
we associate $f$ with $n$.
Conversely, suppose {\rm (ii)} holds, and fix a divisor $n$ of $(q+1)/4$.  
Then there is a bijection between
\begin{itemize}
\item the set of equivalence classes of polynomials $f$ which satisfy
{\rm (i)} and are associated with $n$; and
\item the set of even-order elements in $k^*/(k^*)^{2n}$.
\end{itemize}
One such bijection is defined as follows: for $\alpha\in k^*$, if the coset
$\alpha(k^*)^{2n}\in k^*/(k^*)^{2n}$ has even order, then this coset
corresponds to the equivalence class of the polynomial
\[X(X^{2n}-\alpha)^{(q+1)/(4n)}
\left(\frac{(X^{2n}-\alpha)^{(q-1)/2}+\alpha^{(q-1)/2}}{X^{2n}}\right)^{(q+1)/(2n)}.\]
\end{thm}
In particular, suppose that $k$ is finite and $k\cap\F_q=\F_3$.
Then, for each divisor $n$ of $(q+1)/4$, there is a unique equivalence
class of polynomials $f$ which satisfy (i) and which are associated
with $n$.  For $k=\F_3$, or more generally if $k\nsupseteq\F_9$,
the examples that arise are precisely the polynomials described
in~\cite{LZ}.  If $k\supseteq\F_9$ then the polynomials in
the theorem are new examples of indecomposable exceptional
polynomials.  Finally, let $r$ be the largest
power of 2 dividing $[k\col\F_3]$; then, for any fixed divisor
$n$ of $(q+1)/4$, the equivalence class of polynomials associated
with $n$ and satisfying (i) contains a polynomial defined over
$\F_{3^r}$ but does not contain any polynomials defined over proper
subfields of $\F_{3^r}$.

There exist infinite fields $k$ (for instance, $k=\F_3(y)$ with $y$
transcendental) for which there are infinitely many equivalence
classes of polynomials over~$k$ satisfying (i).  It is interesting
to note, however, that for any $k$ (finite or infinite), any
polynomial over $k$ which satisfies (i) is equivalent over $\bar{k}$
to one of the polynomials over $\F_3$ exhibited in~\cite{LZ},
even though the latter polynomial might not be exceptional over~$k$.
We will see below that a similar remark applies in the case of
characteristic 2.

We now prove Theorem~\ref{char3thm}.  We first show that exceptionality
cannot hold if 
$G:=\Gal(f(X)-u,\bar{k}(u))\cong\PGL(q)$: assume the opposite.
Arguing as in the proof of Lemma~\ref{zero}, we can identify
$A:=\Gal(f(X)-u,k(u))$ with a subgroup of $\PGammaL(q)$,
so $A/G$ is cyclic.  Exceptionality of $f$ implies:
for any $\mu\in A$, if $\mu G$ generates $A/G$ then $\mu$ has exactly one
fixed point (cf.\ \cite[Lemma~6]{Coh} or \cite[Lemma 4.3]{GTZ}).  Since $q\equiv 3\pmod{4}$,
we can view the permutation actions
of $G$ and $A$ as actions on the involutions in $L:=\PSL(q)$ (as in the
proof of Lemma~\ref{fps}).  For any $\mu\in A$ such that $\mu G$ 
generates $A/G$, note that $\mu G=\sigma G$ for some field automorphism
$\sigma$, and $\sigma$ centralizes $\PSL(3)$ and so fixes at least
three involutions of $L$, contradicting exceptionality.

In light of Theorem~\ref{groups}
it follows that any indecomposable exceptional polynomial as in
Theorem~\ref{char3thm} must be equivalent over $\bar{k}$ to one
of the polynomials in Theorem~\ref{genpol2} with $m$ even.
There are several ways to
proceed; the quickest, given what has already been done, is 
to compute the fields of definition for the irreducible factors of
$f(X)-f(Y)$ in $\bar{k}[X,Y]$, for each polynomial $f$ over
$\bar{k}$ which is equivalent to one of the polynomials in
Theorem~\ref{genpol2} with $m=2n$.
This bivariate factorization was derived in~\cite{Z}.
This gives the statement of Theorem~\ref{char3thm}.
The result (and its proof) for characteristic two is similar.  
\begin{thm}
\label{char2thm}
Let $k$ be a field of characteristic $2$, and let $q=2^e$
with $e>1$ odd.  The following are equivalent:
\begin{enumerate}
\item there exists an indecomposable exceptional polynomial $f\in
k[X]$ of degree $q(q-1)/2$ for which\/ $\PGL(q)$ is a transitive normal
subgroup of $\Gal(f(X)-u,k(u))$, and the extension
$k(x)/k(f(x))$ either has no finite branch points or has Galois closure
of genus $\ne (q^2-q)/2$;
\item $k\cap\F_q=\F_2$.
\end{enumerate}
Moreover, for any $f$ as in {\rm (i)}, the Galois closure of
$\bar{k}(x)/\bar{k}(f(x))$ is $\bar{k}(v,w)$ where $v^q-v=w^n$ and $n$
divides $q+1$;
here $n$ is uniquely determined by $f$, and we associate $f$ with $n$.
Conversely, suppose {\rm (ii)} holds, and fix a divisor $n$ of $q+1$.
Then there is a bijection between%
\footnote{This statement must be modified
in case $n=q+1$ and $k$ is not perfect.  In this case, we will define a
bijection between the set of $k$-equivalence classes of
polynomials $f$ as in {\rm (i)} which are associated with $q+1$,
and the set of equivalence classes of pairs $(\alpha,\beta)\in k^*\times k$
modulo the relation: $(\alpha,\beta)\sim (\alpha',\beta')$ if there exists $\gamma\in k$
such that $\gamma^{q+1}=\alpha'/\alpha$ and $\sqrt{\beta}+\gamma\sqrt{\beta'}\in k$. Specifically,
we let the class of $(\alpha,\beta)$ correspond to the $k$-equivalence class of
the polynomial
\[\beta+\sum_{i=0}^{e-1} \alpha^{2^i-1}(X+\beta)^{q(2^i-1)+2^i}.\]}
\begin{itemize}
\item the set of equivalence classes of polynomials $f$ which satisfy
{\rm (i)} and are associated with $n$; and
\item $k^*/(k^*)^n$.
\end{itemize}
One such bijection is defined as follows: for $\alpha\in k^*$, the coset
$\alpha(k^*)^n\in k^*/(k^*)^n$ corresponds to the equivalence class of
the polynomial
\[X\biggl(\sum_{i=0}^{e-1}(\alpha X^n)^{2^i-1}\biggr)^{(q+1)/n}.\]
\end{thm}
Consider the case of finite $k$, and assume $k\cap\F_q=\F_2$.
Then each divisor $n$ of $q+1$ corresponds to a unique
$k$-equivalence class of polynomials as in (i), except when
$k\supseteq\F_4$ and $3\mid n$, in which case there are three
classes.  For $k=\F_2$, or more generally if either $k\nsupseteq\F_4$
or $3\nmid n$, this implies that the only polynomials satisfying (i)
are the ones described in~\cite{CM,Mu}.  However, if
$k\supseteq\F_4$ and $3\mid n$ then we get new examples of
indecomposable exceptional polynomials.  

Finally, consider Theorem~\ref{ind}.  The proofs in \cite{GS,GS2}
show that, in cases (ii) and (iii), the arithmetic monodromy group
$A:=\Gal(f(X)-u,k(u))$ is either $\PSL(p)$ or $\PGL(p)$,
with $p\in\{7,11\}$.  Indecomposability of $f$ implies that $A$
is a primitive permutation group of degree $(p^2-p)/2$.
Since $f$ decomposes over some extension of~$k$, it
decomposes over~$\bar k$, so the geometric monodromy
group $G:=\Gal(f(X)-u,\bar{k}(u))$ is a normal imprimitive
subgroup of $A$.  Thus we must have $A=\PGL(p)$ and $G=\PSL(p)$.
Now Theorems~\ref{groups} and \ref{genpol2} determine $f$ up to
equivalence over $\bar{k}$, and a straightforward computation
with coefficients determines which members of these equivalence
classes are indecomposable over $k$.  The result is as follows:

\begin{thm}
\label{indthm}
Let $k$ be a field of characteristic $p\ge 0$.
The following are equivalent:
\begin{enumerate}
\item There exists an indecomposable polynomial $f\in k[X]$ 
such that $f$ decomposes over some extension of $k$ and
$\deg(f)$ is not a power of $p$;
\item $k$ contains nonsquares and $p\in\{7,11\}$.
\end{enumerate}
Moreover, for any $f$ as in~{\rm (i)}, the Galois closure of 
$\bar{k}(x)/\bar{k}(f(x))$ is $\bar{k}(v,w)$ where $v^p-v=w^n$ 
and $n\mid 2$; here $n$ is uniquely determined by $f$, and
we associate $f$ with $n$.
Conversely, suppose {\rm (ii)} holds, and fix a divisor $n$ of $2$.
Then there is a bijection between
\begin{itemize}
\item the set of equivalence classes of polynomials $f$ which satisfy
{\rm (i)} and are associated with $n$; and
\item the set of nonsquares in $k^*/(k^*)^{2n}$.
\end{itemize}
One such bijection is defined as follows: for any nonsquare $\alpha\in k^*$,
the coset $\alpha(k^*)^{2n}\in k^*/(k^*)^{2n}$ corresponds to
the $k$-equivalence class of
\[{}\begin{cases}
X(X^{2n}-\alpha)^{2/n}\left(\frac{(X^{2n}-\alpha)^3+
\alpha^3}{X^{2n}}\right)^{4/n}&\text{if $p=7$}\\
X^3 (X^{2n}+3\alpha)^{2/n} h(X^{2n})^{6/n}&\text{if $p=11$},
\end{cases}
\]
where $h(X) = X^4 - \alpha X^3 + 2 \alpha^2 X^2 - \alpha^3 X - 3 \alpha^4$.
\end{thm}
%
%

\subsection{Proof of Theorem~\ref{genpol2}}
\label{subsec}

Let $f(X)\in k[X]$ satisfy the hypotheses of Theorem~\ref{genpol2}, 
let $q=p^e$, and assume that item (1) of the theorem holds.  
Then Theorem~\ref{groups} applies.  Since the conclusion of
Theorem~\ref{groups} is unchanged when we replace $k$ by $\bar{k}$,
and likewise this replacement does not affect (2) of Theorem~\ref{genpol2},
we assume henceforth that
\[\textit{$k$ is algebraically closed}.\]

Let $E$ denote the Galois closure
of $k(x)/k(f(x))$, let $G=\Gal(E/k(f(x))$, and let
$G_1=\Gal(E/k(x))$.  Let $I$ be the inertia group of a place $P$
of $E$ lying over the infinite place of $k(f(x))$, and let
$V$ be the Sylow $p$-subgroup of $I$.
Assume that the triple $q,G,G_1$ is not listed in Table~A (we will
return to the cases in Table~A at the end of the proof).
Then $f$ satisfies (i)-(iv) of Theorem~\ref{groups}.
In particular, $G$ is either $\PGL(q)$ or $L$, and $I$ is a Borel
subgroup of $G$, and the higher ramification groups for $P$ satisfy
$V=I_1=\dots=I_n\gneqq G_{n+1}=1$.  Also we have $q>3$.
We first show that $E=k(v,w)$ where $v^q-v=w^n$.

\begin{lemma} There exists $t\in E$ such that
$E^I=k(t)$ and the ramification in $E/E^I$ is as
in {\em($\dagger$)} of Lemma~\ref{lem}, with $r=q-1$ if $G=\PGL(q)$
and $r=(q-1)/2$ otherwise.
\end{lemma}

\begin{proof}
First we show that $I$ has trivial intersection with the inertia
group $J$ of any place of $E$ lying over a finite place
of $k(f(x))$.  By (iv) of
Theorem~\ref{groups}, any such $J$ is cyclic of order
dividing $q+1$; since $\gcd(q+1,|I|)\leq 2$, if $J\cap I\neq 1$ then 
$J\cap I$ has order $2$.  If $q\equiv 3\!\pmod{4}$ then 
$2\mid |I|$ implies that $G=\PGL(q)$, and all involutions in $I$ 
lie in $G\setminus L$;
but the order of $J$ is divisible by 4, so its
involution lies in $L$, whence $J\cap I=1$.  If $q\equiv1\!\pmod{4}$,
then the involution of $J$ lies in $G\setminus L$ but the
involutions of $I$ lie in $L$.  If $q\equiv 0\!\pmod{2}$, then
$J$ contains no involution.   Thus, in all cases $J\cap I=1$,
so $E/E^I$ is unramified over places of $E^I$
lying over finite places of $k(f(x))$.

Next we consider the inertia groups of places of $E$ lying over
the infinite place of $k(f(x))$; these are precisely the
conjugates of $I$ in $G$ (i.e.\ the Borel subgroups of $G$), and, since $I$ 
is self-normalizing in $G$, no two of these places have the same 
inertia group.  The intersection of any two distinct Borel subgroups of $G$
has order $r$ (by Lemma~\ref{borels}), so each of the $q$ places
of $E$ which differ from $P$ and lie over the infinite place
of $k(f(x))$ have ramification index $r$ in
$E/E^I$.  Hence all ramification in $E/E^I$
occurs over two places of $E^I$, one of which is totally
ramified and the other of which has ramification index $r$.
We now compute the genus of $E^I$:
\[(q-1)(n-1)-2=2g_{E}-2=(2g_{E^I}-2)qr+(qr-1)+n(q-1)+q(r-1),\]
so $g_{E^I}=0$.  Since $k$ is algebraically closed, it
follows that $E^I=k(t)$ for some $t$.
Finally, we can replace $t$ by some $(\alpha t + \beta)/(\gamma t + \delta)$
(with $\alpha,\beta,\gamma,\delta\in k$) to make the infinite place of $k(t)$ be totally ramified
in $E$, and to make $0$ be the only finite place of $k(t)$ which ramifies in $E/k(t)$.
\end{proof}

For $n\neq q+1$, we will apply Theorem~\ref{unq} to the extension 
$E/E^I$; 
first we must verify (i) and (ii) of that result.  
Condition (ii) is trivial, since $r/\gcd(n,r)\geq (q-1)/2>\sqrt{q}-1$.
For the specific values
of $n$ and $r$ under consideration, condition (i) is proved
in Lemma~\ref{numerology}.  Hence Theorem~\ref{unq} applies to
the extension $E/E^I$, so $E=k(v,w)$ where
$v^q-v=w^n$.  

Now assume $n=q+1$, so $p=2$ and $g_{E}=(q^2-q)/2$.  
The previous lemma shows that $E/E^V$ is only ramified
over one place, where it is totally ramified with the only
jump in the ramification occurring after the $n$-th ramification
group; this implies that $E^V$ has genus zero,
so $E^V=k(\tilde w)$.  Then Proposition~\ref{msri} applies to
the tower $E\supseteq k(\tilde w)\supseteq E^V$, so we
conclude that $E=k(v,w)$ where $v^q-v\in k[w]$; the degree
$n'$ of any term of $v^q-v$ (as a polynomial in $w$)
satisfies $n'\equiv q+1\!\pmod{q-1}$ and $q\nmid n'$, and moreover
$q+1$ is the largest integer coprime to $p$ which divides some such
$n'$.  We can write $n'=r 2^i$ where $r$ is odd,
$1\leq r\leq q+1$, and $0\leq i<e$.  Then we have
$r\equiv 2^{e+1-i}\!\pmod{2^e-1}$, and the only possibilities
are $r=2^e+1$ (and $i=0$) or $r=1$ (and $i=1$).  Thus,
$v^q-v=\alpha w^2+\beta w^{q+1}$, where $\alpha\in k$ and $\beta\in
k^*$ (and, by replacing $w$ by $\beta^{1/(q+1)}w$, we may assume 
$\beta=1$).
If $\alpha\neq 0$ then, by~\cite[Satz 7]{St2}, every $k$-automorphism
of $E$ preserves the place $P$; hence $\Aut_k(E)$ is
the decomposition group of $P$ in the extension
$E/E^{\Aut_k(E)}$, so it is solvable, and thus has
no subgroup isomorphic to $L$, a contradiction.  Thus
$E=k(v,w)$ where $v^q-v=w^{q+1}$.

Our next task is to determine all subgroups of 
${\mathcal G}:=\Aut_k(E)$ isomorphic to $G$.  
Recall that $n$ divides $q+1$.  We exhibit three subgroups of
$\mathcal G$.  There is an elementary abelian subgroup $U$ of
order $q$, whose elements fix $w$ and map $v\mapsto v+\alpha$
with $\alpha\in\F_q$.  There is a cyclic subgroup $H$ of order
$n(q-1)$, whose elements map $v\mapsto\zeta^n v$ and
$w\mapsto\zeta w$ where $\zeta^{n(q-1)}=1$.  And there is a cyclic subgroup $J$
of order $2$ or $4$, generated by the automorphism sending
$v\mapsto 1/v$ and $w\mapsto (-1)^{1/n}w/v^{(q+1)/n}$,
for any choice of $n$-th root of $-1$; this group has order $2$
precisely when $(-1)^{1/n}=-1$.  One can verify that these maps
are automorphisms by observing that they are bijections of
$k(v)[X]$ which induce bijections on the ideal generated by
$v^q-v-X^n$, so they are bijective on the quotient $E$.
Let ${\mathcal G}_0$ be the group generated by $U$, $H$, and $J$.
Note that each element of $U\cup H\cup J$ induces an
automorphism of $k(v)$; the induced automorphisms are
$v\mapsto v+\alpha$ (for $\alpha\in\F_q$), $v\mapsto\hat\zeta v$
(for $\hat\zeta\in\F_q^*$), and $v\mapsto 1/v$.  These automorphisms
generate the group $\PGL(q)$.  Hence, restriction to $k(v)$
induces a homomorphism $\rho\colon{\mathcal G}_0\to\Aut_k(k(v))$
whose image is $\PGL(q)$; the kernel of $\rho$ is
the cyclic subgroup $Z$ of $H$ of order $n$ (since
$Z=\Gal(E/k(v))$), so $|{\mathcal G}_0|=n(q^3-q)$.  
Also, $Z$ commutes with all 
elements of $U$, $H$, and $J$, so $Z$ lies in the center of 
${\mathcal G}_0$.

Now we compute the subfield $E^{{\mathcal G}_0}$; this field
equals $k(v)^{\PGL(q)}$.  The latter field was
computed by Dickson~\cite[p.~4]{Di2}: it is $k((v^{q^2}-y)^{q+1}/
(v^q-v)^{q^2+1})$ (to verify this, note that this rational
function is fixed by each element of $\PGL(q)$, and its degree
is $q^3-q=|\PGL(q)|$).  Since $v^q-v=w^n$, we can rewrite
this generator in terms of $w$: 
$E^{{\mathcal G}_0}=k(w^{n(q-q^2)}(w^{n(q-1)}+1)^{q+1})$.
Put $\hat u:=w^{q-q^2}(w^{n(q-1)}+1)^{(q+1)/n}$; then we have
$E^{{\mathcal G}_0}=k({\hat u}^n)$, so $\hat{G}:=\Gal(E/k(\hat u))$
is a subgroup of ${\mathcal G}_0$ of order $q^3-q$.

Suppose $p$ is odd.
Put $\tilde u:=w^{(q-q^2)/2}(x^{n(q-1)}+1)^{(q+1)/(2n)}$, so
$\hat L:=\Gal(E/k(\tilde u))$ is a subgroup of $\hat G$ of index two.
Then $\tilde u$ is not fixed by any nontrivial element of $Z$, so
$\rho$ induces an isomorphism between $\hat L$ and the unique subgroup
of $\PGL(q)$ of index two (namely $\PSL(q)$).  Moreover, $\hat L$ is
a normal subgroup of ${\mathcal G}_0$ with cyclic quotient of
order $2n$.

Now suppose $n$ is odd (but we no longer restrict $p$).  Then
$\hat u$ is not fixed by any nontrivial element of $Z$,
so $\rho$ induces an isomorphism $\hat G\cong\PGL(q)$; thus
${{\mathcal G}_0}=\hat G\times Z\cong\PGL(q)\times C_n$.
Note that all elements of $\hat G\times Z$ of order
$p$ are in $\hat G$; thus, any
subgroup of ${\mathcal G}_0$ isomorphic to $\PSL(q)$ contains all
such elements, and so must be the unique subgroup of $\hat G$ isomorphic
to $\PSL(q)$ (since this subgroup is generated by the elements of
order $p$).  If $p=2$ then this subgroup is $\hat G$.  If $p>2$ then this
subgroup is $\hat L$.  Moreover, if $1<n<q+1$ and $p>2$, then ${\mathcal G}_0/{\hat L}$
is cyclic of order $2n$, so $\hat G$ is the only subgroup of ${\mathcal G}_0$
isomorphic to $\PGL(q)$.

Next assume that $n$ is even.  Then $n$ divides
$(q+1)/\gcd(4,q+1)$, so $q\equiv 7\!\pmod{8}$.
In particular, $q$ is odd, so $1<n<q+1$.  From above we know that
${\mathcal G}_0$ has
a unique subgroup isomorphic to $\PSL(q)$.  We now show that 
${\mathcal G}_0$ has no subgroup isomorphic to $\PGL(q)$.  We know that ${\mathcal G}_0$
has a normal subgroup $\hat L$
with cyclic quotient of order $2n$.  
If ${\mathcal G}_0$ has a subgroup $\tilde G$ isomorphic to $\PGL(q)$,
then (since $\PGL(q)$ has trivial center) we would have
${\mathcal G}_0={\tilde G}\times Z$, so ${\mathcal G}_0/{\hat L}\cong C_2\times C_m
\not\cong C_{2m}$, a contradiction.

If $1<n<q+1$, then Stichtenoth proved that
$|{\mathcal G}|=n(q^3-q)$~\cite{St2}, so ${\mathcal G}={\mathcal
G}_0$; hence $\mathcal G$ has a unique subgroup isomorphic to
$\PSL(q)$, and has a subgroup isomorphic to $\PGL(q)$ if and only if
$n$ is odd, in which case it has a unique such subgroup.
If $n=1$ then $E=k(v)$, so ${\mathcal G}=\PGL(k)$; this group
has a unique conjugacy class of subgroups isomorphic
to $\PSL(q)$, and a unique conjugacy class of subgroups isomorphic
to $\PGL(q)$.  Finally, if $n=q+1$ (so $p=2$) then Leopoldt showed
that ${\mathcal G}=\PGU(q^2)$~\cite{Le}, and this group has a 
unique conjugacy class of subgroups isomorphic to $\PGL(q)$.
Since conjugate groups will lead to the same polynomials, we may assume
in every case that
\[E^G\, \text{ is either $k(\hat u)$ or $k(\tilde u)$}.\]

We now compute the subfield of $E$ invariant under the
one-point stabilizer $G_1$ of $G$.  Since $G$ has a unique conjugacy
class of subgroups isomorphic to $G_1$, and conjugate groups $G_1$ will
lead to the same polynomials, it suffices to do this for a single
point-stabilizer $G_1$.  

To complete the proof, we must compute the polynomials $f$ in each 
of three cases:
$q$ is odd and $G=\hat G\cong\PGL(q)$; $q\equiv3\!\pmod{4}$
and $G=\hat L\cong\PSL(q)$; $q$ even and $G=\hat G\cong\PGL(q)$.
The computations in each case are similar, so we only give the
details for the first case.  Thus, for the remainder of the proof
we assume that 
\[\text{$q$ is odd and $G=\hat G\cong\PGL(q)$,}\]
so also
\[\text{$n$ is odd and $G_1$ is dihedral of
order $2(q+1)$.}\]
To exhibit such a group $G_1$, view $\GL(q)$ 
as the invertible $\F_q$-linear
maps on a two-dimensional $\F_q$-vector space, and choose this
vector space to be $\F_{q^2}$.  Then the multiplication maps by the
various elements of $\F_{q^2}^*$ form a cyclic subgroup of $\GL(q)$;
the group generated by this cyclic group and the $q$-th power map has
order $2(q^2-1)$, and its image in $\PGL(q)$ is dihedral of 
order $2(q+1)$.  We make this explicit by choosing a nonsquare
$\gamma_0\in\F_q^*$ and letting $\delta\in\F_{q^2}^*$ be a square root
of $\gamma_0$.  Choose the basis $\{1,\delta\}$ for $\F_{q^2}/\F_q$;
with respect to this basis, the matrix for the multiplication map
by $\alpha+\beta\delta$ (where $\alpha,\beta\in\F_q$ are not both zero) is
$\textmatrix{\alpha&\beta\gamma_0\\ \beta&\alpha}$, and the matrix for the $q$-th
power map is $\textmatrix{1&0\\0&-1}$.  These matrices generate
a dihedral subgroup of $\PGL(q)$ of order $2(q+1)$; let $G_1$ be the
intersection of $\hat G$ with
the preimage under $\rho$ of this dihedral group.  

Our next task is to compute the subfield of $k(v)$ fixed by $G_1$;
this will coincide with $E^{G_1\times Z}$.
Let $\gamma=1/\gamma_0$.  Certainly $G_1$
fixes the sum of the images of $v^2$ under $G_1$, namely
\[
2\hat x:=
\frac{1}{q-1}\sum_{\substack{\alpha,\beta\in\F_q\\ (\alpha,\beta)\neq(0,0)}}
 \left(\left(\frac{\alpha v + \beta}{\beta \gamma_0 v + \alpha}\right)^2+
 \left(\frac{\alpha v - \beta}{\beta \gamma_0 v - \alpha}\right)^2\right).
\]
Here the factor $1/(q-1)$ comes from the fact that the pairs
$(\alpha,\beta)$ and $(\alpha\zeta,\beta\zeta)$ (with $\zeta\in\F_q^*$) correspond to the
same element of $G_1$.  Thus each term of $G_1$ corresponds to
a unique pair $(\alpha,\beta)\in\F_q^2$ with either $\beta=\gamma$ or
both $\beta=0$ and $\alpha=1$, so
\begin{equation}
\label{ohy}
\hat x =
v^2+\sum_{\alpha\in\F_q}\left(\frac{{\alpha v + \gamma}}{v + \alpha}\right)^2;
\end{equation}
note that $\hat x=b(v)/(v^q-v)^2$ for some monic $b(X)\in\F_q[X]$
of degree $2q+2$.  Here $b(X)$ has no term of degree $2q+1$.  
For any $\alpha\in\F_q$, multiply~(\ref{ohy}) by $(v+\alpha)^2$ and
then substitute $v=-\alpha$ to find $(-\alpha^2+\gamma)^2=b(-\alpha)$;
hence $b(X)=(X^2-\gamma)^2+(X^q-X)c(X)$ where $c(X)\in\F_q[X]$
is a monic polynomial of degree $q+2$ having no term of degree $q+1$.
We compute the derivative: $b'(X)=4X(X^2-\gamma)-c(X)+
(X^q-X)c'(X)$.
For any $\alpha\in\F_q$, multiply~(\ref{ohy}) by $(v+\alpha)^2$, take the
derivative of both sides, and then substitute $v=-\alpha$; 
this gives $2\alpha(-\alpha^2+\gamma)=b'(-\alpha)$, so
$c(v)=2v(v^2-\gamma)+(v^q-v)a(v)$ where $a(X)=X^2+\xi$ and $\xi\in\F_q$.
The choice of $\xi$ is irrelevant, since changing $\xi$ amounts to adding
a constant to $\hat x$; but regardless,
it is easy to show that $\xi=0$ (e.g.\ by computing appropriate
terms of $(v^q-v)^2{\hat x}$).  Hence 
\[b(v)=\left((v^2-\gamma)+v(v^q-v)\right)^2=
\left(v^{q+1}-\gamma\right)^2,\]
so 
$\hat x={(v^{q+1}-\gamma)^2/(v^q-v)^2}$.
Since $[k(v)\col k(\hat x)]=\deg(\hat x)=2(q+1)$ equals $|G_1|$
and $G_1$ fixes $\hat x$, it follows that $E^{G_1\times Z}=
k(v)^{G_1}=k(\hat x)$.

We have shown that $k(v)^G=k({\hat u}^n)$ and $k(v)^{G_1}=k(\hat x)$,
where ${\hat u}^n=(v^{q^2}-v)^{q+1}/(v^q-v)^{q^2+1}$; hence
${\hat u}^n=h(\hat x)$ for some $h(X)\in k(X)$ which we now 
determine.  Write $h(X)=\xi\prod_i(X-\alpha_i)^{r_i}$,
where the $\alpha_i$ are distinct elements of $k$, 
the $r_i$ are nonnegative integers, and $\xi\in k^*$.
Then we have
\begin{equation}
\label{len}
\frac{\left(\frac{{v^{q^2}-v}}{v^q-v}\right)^{q+1}}
 {(v^q-v)^{q^2-q}}
=\frac{\xi\prod_i\left((v^{q+1}-\gamma)^2-\alpha_i(v^q-v)^2\right)^{r_i}}
  {(v^q-v)^{2\sum_i r_i}}.
\end{equation}
Let $h_i(X)=(X^{q+1}-\gamma)^2-\alpha_i(X^q-X)^2$;
since $\gamma$ is a nonsquare in $\F_q^*$, no two $h_i$
have a common root, and also no $h_i$ has a root in $\F_q$.
Since the poles of the left side of~(\ref{len}) are precisely
the elements of $\F_q$, we conclude that each $r_i$ is positive
(so $h$ is a polynomial) and $\sum r_i=(q^2-q)/2$.  Equating 
the leading coefficients
of the two sides of~(\ref{len}) gives $\xi=1$.  The roots $y=\omega$ of
the left side of~(\ref{len}) are precisely the values
$\omega\in\F_{q^2}\setminus\F_q$, and each has multiplicity $q+1$;
any such $\omega$ is a root of $h_i(X)$ where
\[\alpha_i=\frac{(\omega^{q+1}-\gamma)^2}{(\omega^q-\omega)^2},\]
so $\alpha_i^{(q+1)/2}=-\alpha_i$.  Since $2h_i(X)+(X^q-X)^2h_i'(X)=
2(X^{q+1}-\gamma)(X^2-\gamma)^q$, any multiple root $\eta$ of $h_i$
must satisfy either $\eta^{q+1}=\gamma$ or $\eta^2=\gamma$; in the
former case, $\alpha_i=0$ and $h_i=(X^{q+1}-\gamma)^2$, so $\eta$ is
a root of multiplicity two; in the latter case, $\alpha_i=\gamma$
and $h_i=(X^2-\gamma)^{q+1}$, so $\eta$ is a root of multiplicity
$q+1$.  Hence $h(X)\in k[X]$ divides
\[X^{(q+1)/2}(X-\gamma)
  \left(\frac{X^{(q-1)/2}+1}{X-\gamma}\right)^{q+1};\]
since both these polynomials are monic of degree $(q^2-q)/2$,
they must be the same.

We now determine all possibilities for the original
polynomial $f$.  We have shown that $h(\hat x)={\hat u}^n$,
where $E^{G_1\times Z}=k(\hat x)$ and $E^G=k(\hat u)$;
we will modify this polynomial identity to express $\hat u$
as a polynomial in some $\hat x$ with $E^{G_1}=k(\hat x)$.
We know that $\hat x$ will be the $n$-th root of some 
generator of $k(\hat x)$; from above,
\[\hat x-\gamma=\frac{(v^2-\gamma)^{q+1}}{(v^q-v)^2}
              =\frac{(v^2-\gamma)^{q+1}}{x^{2n}},\]
so we choose $\hat x=(v^2-\gamma)^{(q+1)/n}/w^2$.
Then $\hat h(X):=h(X+\gamma)$ satisfies
$\hat h(\hat x^n)=\hat u^n$; taking $n$-th roots gives
\[\zeta\hat u=
(\hat x^n+\gamma)^{(q+1)/(2n)}{\hat x}
 \left(\frac{(\hat x^n+\gamma)^{(q-1)/2}+1}{\hat x^n}\right)^{(q+1)/n}\]
for some $n$-th root of unity $\zeta$.  We may assume $\zeta=1$
(by replacing $\hat u$ by $\zeta\hat u$).  Thus
$\hat u=\hat f(\hat x)$ where $\hat f(X)\in k[X]$ is defined by
\[\hat f(X):=X(X^n+\gamma)^{(q+1)/2n}
 \left(\frac{(X^n+\gamma)^{(q-1)/2}+1}{X^n}\right)^{(q+1)/n}.\]
Hence $\Gal(E/k(\hat x))$ is an index-$n$ subgroup of $G_1\times
Z$ which is also a subgroup of $G$, so it is $G_1$.  Since $G_1$ 
contains no nontrivial normal subgroups of $G$, the Galois
closure of $k(\hat x)/k(\hat u)$ is $E$, so the monodromy
group of $\hat f$ is $G$.  The only choices we made which restricted
the possibilities for the original polynomial $f$ were the
choices of generators $\hat x$ and $\hat u$ for the fields
$E^{G_1}$ and $E^G$; hence, the polynomials $f$ in this
case are precisely the polynomials $\ell_1\circ {\hat f}\circ\ell_2$,
where $\ell_1$ and $\ell_2$ are linear polynomials
in $k[X]$.
This completes the proof of the theorem in case $q$ is odd and
$G\cong\PGL(q)$.  As noted above, the completion of the proof in
the other two cases is similar.  In the statement of the result there
is only a single family of polynomials covering both cases with $q$ odd;
this occurs because, when $q$ is odd and $G\cong\PSL(q)$, the resulting
polynomials can be obtained by substituting $2n$ for $n$ in the above
expression for $\hat{f}$.

Finally, we consider the group-theoretic possibilities in Table~A,
which lead to the polynomials in Table~B.  The cases with $q$ odd
can be treated in a similar manner to what we have done above; in
particular, in each case a fixed point computation similar to that
of Lemma~\ref{fps} can be used to determine a short list of possibilities
for the ramification
(as in the proof of Theorem~\ref{groups}), after which
Theorem~\ref{unq} implies that the Galois closure of
$k(x)/k(f(x))$ is $k(v,w)$ where $v^q-v=w^n$,
and we conclude the proof precisely as above.  Table~A also includes
the possibility that $q=4$, in which case $G=\PGammaL(q)\cong S_5$.
A Riemann-Hurwitz computation shows that $I$ must be cyclic of order
6, and that the Galois closure $E$ of $k(x)/k(f(x))$ has
genus one.  However, the automorphism groups of function fields of
genus one are known: in our case, $\Aut_k(E)$ is the semidirect
product $H\rtimes J$, where $H$ is the (abelian) group of translations
by points on the elliptic curve corresponding to $E$, and $J$ has
order either $2$ or $24$.  If $\Aut_k(E)$
has a subgroup $G$ isomorphic to $S_5$, then $G\cap H$ is a normal
subgroup of $G$ of index less than 24, so $G\cap H\supseteq A_5$;
in particular, $H$ contains a Sylow $2$-subgroup of $A_5$, namely
a Klein $4$-group, contradicting standard results about the structure
of the group of points on an elliptic curve~\cite[Cor.~6.4]{Si}.
Thus there are no polynomials corresponding to the case $q=4$
in Table~A.  This concludes the proof of Theorem~\ref{genpol2}.

\appendix

\section*{Appendix. Group theoretic preliminaries}

\setcounter{thm}{0}
\renewcommand{\thethm}{A.\arabic{thm}}

In this appendix we summarize the basic group theoretic facts 
used in this paper.
Let $\F_q$ be a field of order $q$ and characteristic $p$.
As usual, $\GL(q)$ denotes the group of invertible
two-by-two matrices over $\F_q$, and $\SL(q)$ is its subgroup
of determinant-one matrices.
The centers of these groups are the scalar
matrices, which in the case of $\GL(q)$ are just $\F_q^*$
and in the case of $\SL(q)$ are $\{\pm1\}$; the quotients of $\GL(q)$
and $\SL(q)$ by their centers are denoted $\PGL(q)$ and $\PSL(q)$,
respectively.  We often write $L$ for $\PSL(q)$.
The orders of these groups are as follows:
$|\GL(q)|=(q^2-1)(q^2-q)$, $|\PGL(q)|=|\SL(q)|=q^3-q$, and
$|\PSL(q)|=(q^3-q)/o$, where $o=\gcd(2,q-1)$.
 
We now discuss the structure of $L=\PSL(q)$.
The most important property of $L$ is that, for $q\geq 4$,
it is a simple group~\cite[\S1.9]{Su}.
For small $q$, we have the isomorphisms $\PSL(2)\cong S_3$,
$\PSL(3)\cong A_4$, and $\PSL(4)\cong\PSL(5)\cong A_5$.
For general $q$, the subgroups of $\PSL(q)$ were 
determined by Dickson in 1901~\cite[\S 260]{Di};
a treatment in modern
language is~\cite[\S3.6]{Su}.  We state the result for the reader's
convenience.  As above, $o=\gcd(2,q-1)$.
\begin{thm}[Dickson, 1901] 
\label{dickson}
The subgroups of $\PSL(q)$ are precisely the following groups.
\begin{enumerate}  
\item The dihedral groups of order $2(q\pm1)/o$ and their
  subgroups.
\item The group $B$ of upper triangular matrices in\/ $\PSL(q)$, and
  subgroups of $B$ (here, the order of $B$ is $q(q-1)/o$, the
  Sylow $p$-subgroup $U$ of $B$ is elementary abelian and normal in
  $B$, and the quotient group $B/U$ is cyclic).
\item $A_4$, except if $q=2^e$ with $e$ odd.
\item $S_4$, if $q\equiv\pm1\pmod 8$.
\item $A_5$, except if $q\equiv\pm2\pmod 5$.
\item $\PSL(r)$, where $r$ is a power of $p$ such that $r^m=q$.
\item $\PGL(r)$, where $r$ is a power of $p$ such that $r^{2m}=q$.
\end{enumerate}
\end{thm}

This result describes the isomorphism classes of subgroups
of $L=\PSL(q)$; we are also interested in conjugacy classes of
subgroups.  We only need this in certain cases.  
\begin{lemma}
\label{dihedrals}
For $q>3$, there are $(q^2-q)/2$ subgroups of $L$ which are dihedral of
order $2(q+1)/o$, and any two of these are conjugate.
Let $J$ be a nontrivial cyclic subgroup 
of $L$ whose order divides $(q+1)/o$; then the
normalizer of $J$ in $L$ is dihedral of order $2(q+1)/o$.
\end{lemma}
\begin{proof}
By Theorem~\ref{dickson}, $L$ has a dihedral subgroup $H$ of
order $2(q+1)/o$, and this subgroup is maximal unless $q=7$
(when $S_4$ is a possible overgroup) or $q=9$ (when $A_5$ is a
possible overgroup).  Since $L$ is simple, it follows that $H$
is the normalizer (in $L$) of any of its cyclic subgroups of
order more than two, and also that $H$ is self-normalizing in $L$.
Let $H'$ be another dihedral subgroup
of $L$ of order $2(q+1)/o$, and suppose $|J|>2$.
Let $s$ be an odd prime divisor of $q+1$ if any such exists,
and otherwise put $s=2$; then Sylow $s$-subgroups of 
$H$ and $H'$ are also Sylow subgroups of $L$, and so are conjugate,
whence their normalizers $H$ and $H'$ are conjugate as well.  
Likewise, a Sylow $s$-subgroup of $J$ is conjugate to a cyclic subgroup
of $H$, so its normalizer is conjugate to $H$; since $J$ is contained
in the normalizer, $J$ is conjugate to a subgroup of $H$, so its
normalizer is conjugate to $H$.
We have shown that there is a unique conjugacy
class of dihedral subgroups of $L$ of order $2(q+1)/o$,
so the size of this class is $|L\col N_L(H)|=|L\col H|=(q^2-q)/2$.  
This concludes the proof in case $|J|>2$; the case
$|J|=2$ is included in the next lemma.
\end{proof}
By an \emph{involution} of a group, we mean an element of order two.

\begin{lemma}
\label{involutions}
$L$ contains a unique conjugacy class of involutions.  For
$q$ odd,\/ $\PGL(q)$ contains two conjugacy classes of involutions.
The number of involutions in\/ $\PGL(q)$ is either $q^2$ (if $q$ odd)
or $q^2-1$ (if $q$ even).  For $q$ odd, the number of involutions in
$L$ is either $(q^2-q)/2$ (if $q\equiv3\pmod{4}$) or
$(q^2+q)/2$ (if $q\equiv1\pmod{4}$).  For $q$ even, the centralizer
of an involution of\/ $\PGL(q)$ is a Sylow $2$-subgroup.  For $q$ odd, the
centralizer in\/ $\PGL(q)$ of an involution of\/ $\PGL(q)$
is dihedral of order either $2(q+1)$ or $2(q-1)$, and the centralizer
in $L$ of this involution is dihedral of half the size; order
$2(q+1)$ occurs when $q\equiv 3\pmod{4}$ and the 
involution is in $L$, and also when $q\equiv1\pmod{4}$ and the
involution is not in $L$, and order $2(q-1)$ occurs otherwise.
\end{lemma}
\begin{proof}
One immediately verifies that a nonidentity element 
$\textmatrix{\alpha & \beta \\ \gamma & \delta}$ of $\PGL(q)$ is an involution if and only
if $\delta = -\alpha$; it is then a triviality to count the involutions in
either $\PGL(q)$ or $L$.  For $q$ even, let $U$ be the group
of matrices $\textmatrix{1&*\\0&1}$, which is a Sylow $2$-subgroup
of $L$; then $U$ is elementary abelian and is the centralizer (in $\PGL(q)$) of
any nontrivial $\mu\in U$.  Since any involution in $L$ is
conjugate to an involution in $U$, the centralizer of an involution
is a Sylow $2$-subgroup.  But the number of involutions of $L$ which are 
conjugate to a fixed involution $\mu\in U$ is $|L\col N_L(\langle \mu\rangle)|=|L\col U|=q^2-1$, so 
$L$ contains a unique conjugacy class of involutions.  Henceforth
assume $q$ odd.  First consider the involution
$\textmatrix{-1&0\\0&1}$, which lies in $L$ if and only if 
$q\equiv1\!\pmod4$: its centralizer in $\PGL(q)$ (respectively, $L$) 
consists of all the diagonal and antidiagonal elements, and so is dihedral
of order $2(q-1)$ (respectively, $q-1$).  Hence the number of conjugates
of this involution by either $L$ or $\PGL(q)$ is
$|L|/(q-1)=(q^2+q)/2$.
We will show that $\PGL(q)$ contains a dihedral 
subgroup $H$ of order $2(q+1)$; we now show how the remainder of the 
lemma follows from this statement.  Since $L$ has no cyclic subgroup
of order $q+1$, the intersection $H\cap L$ is dihedral of order $q+1$.
The center of $H$ is generated by an involution $\nu$, and
$\nu$ lies in $L$ if and only if $4\mid(q+1)$.
The centralizer of $\nu$ in $L$ 
contains $H\cap L$; this centralizer is a proper
subgroup of $L$ (since no nonidentity element of $\PGL(q)$
centralizes $L$), so it must be $H\cap L$ (since $H\cap L$ is maximal
unless $q=7$ or 9, and in those cases we note that $S_4$ and $A_5$
have trivial center).  Then the centralizer of $\nu$ in $\PGL(q)$
must be $H$ (since it contains $H$ and has $H\cap L$ as a subgroup of
index at most two).  Finally, the number of conjugates of $\nu$ by
either $L$ or $\PGL(q)$ is $|L\col H\cap L|=|\PGL(q)\col H|=(q^2-q)/2$,
which completes the proof.

Now we must show that $\PGL(q)$ contains a dihedral subgroup of order
$2(q+1)$.  For this, we view $\GL(q)$ as the invertible $\F_q$-linear
maps on a two-dimensional $\F_q$-vector space, and we choose this
vector space to be $\F_{q^2}$.  Then the multiplication maps by the
various elements of $\F_{q^2}^*$ form a cyclic subgroup of $\GL(q)$;
the group generated by this cyclic group and the $q$-th power map has
order $2(q^2-1)$, and its image in $\PGL(q)$ is dihedral of 
order $2(q+1)$.
\end{proof}

Let $G$ denote either $\PGL(q)$ or $\PSL(q)$; then a \emph{Borel subgroup}
of $G$ is any subgroup conjugate to the group of upper-triangular matrices
in $G$.
\begin{lemma}
\label{borels}
Let $G$ be either\/ $\PGL(q)$ or\/ $\PSL(q)$.  A Borel subgroup of $G$
is self-normalizing.  The intersection of any two distinct Borel 
subgroups of $G$ is cyclic of order either $q-1$
(if $G=\PGL(q)$) or $(q-1)/2$ (if $G\neq\PGL(q)$).
\end{lemma}
\begin{proof}
One easily verifies that the upper-triangular matrices $B$ are
self-normalizing, so the same is true of any Borel subgroup of $G$.  Hence there
are $|G\col B|=q+1$ Borel subgroups of $G$.  One of these is the 
lower-triangular matrices, which intersects $B$ in the diagonal.
The number of conjugates of the lower-triangular matrices under $B$
is then the index of the diagonal in $B$, namely $q$, so the
conjugation action of $G$ on its Borel subgroups is doubly
transitive, and thus the intersection of any two distinct Borel subgroups is
a conjugate of the diagonal.
\end{proof}

Next we consider the automorphisms of $\PSL(q)$.  Since $\PSL(q)$
is a normal subgroup of $\PGL(q)$, these include conjugation
by any element of $\PGL(q)$; since no nonidentity element of
$\PGL(q)$ centralizes $\PSL(q)$, this gives an embedding of
$\PGL(q)$ into the automorphism group $\Aut(\PSL(q))$.
Any $\sigma\in\Gal(\F_q/\F_p)$ induces an automorphism of
$\PSL(q)$ by acting on the matrix entries; these automorphisms are
called \emph{field automorphisms}.  The group
$\Gal(\F_q/\F_p)$ normalizes $\PGL(q)$; let $\PGammaL(q)$
denote the semidirect product of $\PGL(q)$ by $\Gal(\F_q/\F_p)$.
Then $\PGammaL(q)=\Aut(\PSL(q))$~\cite[Thm.~12.5.1]{Ca}.
 
We need one result about subgroups of $\PGammaL(q)$.
\begin{lemma}  
\label{unique_action}
Let $G$ be a subgroup of\/ $\PGammaL(q)$ containing $L$, where $q>3$.
There are precisely $(q^2-q)/2$ subgroups $G_1$ of
$G$ such that $|G\col G_1|=(q^2-q)/2$ and $G_1\cap L$ is dihedral of order 
$2(q+1)/o$; any two such subgroups $G_1$ are conjugate.
\end{lemma}

\begin{proof}
Let $\Lambda$ denote the set of dihedral subgroups of $L$ of order
$2(q+1)/o$.  By Lemma~\ref{dihedrals}, we have $|\Lambda|=q(q-1)/2$,
and moreover $G$ acts transitively on $\Lambda$ by conjugation.
Thus, for any $H\in\Lambda$, we have $(q^2-q)/2=|\Lambda|=|G\col N_G(H)|$.
Since any two groups in $\Lambda$ are conjugate, also their
normalizers in $G$ are conjugate.  Conversely, if the subgroup $G_1$ of $G$
satisfies $|G\col G_1|=(q^2-q)/2$ and $G_1\cap L\in\Lambda$, then $G_1$ normalizes
$H:=G_1\cap L$, so we must have $G_1=N_G(H)$.  The result follows.
\end{proof}

The main result on conjugacy classes of $\PGammaL(q)$ is the 
following consequence of Lang's theorem on algebraic 
groups~\cite[2.7-2]{GL}. 
 Let $\sigma\in\PGammaL(q)$ be a field
automorphism and $\mu\in\PGL(q)$; if $\mu\sigma$ and
$\sigma$ have the same order, then
$\mu\sigma=\tau\sigma\tau^{-1}$ for some
$\tau\in\PGL(q)$.  We refer to this result as Lang's theorem.
In the proof of Lemma~\ref{grz}, we need the following more
general consequence of Lang's theorem: if $\sigma$ is a field
automorphism in $\PGammaL(\bar{k})$, then for every $\mu\in\PGL(\bar{k})$
there exists $\tau\in\PGL(\bar{k})$ such that $\tau(\mu\sigma)\tau^{-1}=\sigma$.

We also use Zsigmondy's theorem on primitive
prime divisors.  For a prime $p$ and positive integer $e$, we
say that a prime $s$ is a primitive prime divisor of $p^e-1$
if $e$ is the least positive integer $i$ for which $s$ divides
$p^i-1$.  Zsigmondy's theorem says that, for fixed $p$ and $e$,
there exists a primitive prime divisor of $p^e-1$ unless either
$e=2$ (and $p+1$ is a power of 2) or $p=2$ (and $e=6$).
For a proof, see~\cite[Thm.~6.2]{Lu} 
or the original~\cite{Zs}.  Note that a primitive prime divisor
of $p^e-1$ is coprime to $e$.

\end{document}